\documentstyle[12pt]{article}

\begin{titlepage}

\title{Cohomological properties of the smooth globalization of a Harish-Chandra
module}

\author{Ulrich Bunke\thanks{Humboldt-Universit\"at zu Berlin, Institut f\"ur
Reine Mathematik (SFB288), Ziegelstr. 13a, Berlin 10099.
E-mail:ubunke@mathematik.hu-berlin.de
} and
Martin
Olbrich\thanks{Humboldt-Universit\"at zu Berlin, Institut f\"ur Reine
Mathematik (SFB288), Ziegelstr. 13a, Berlin 10099.
E-mail:olbrich@mathematik.hu-berlin.de  }
}

\end{titlepage}

\newcommand{\proof}{{\it Proof.$\:\:$}}

\newcommand{\dist}{{\rm dist}}
\newcommand{\kaaa}{{\bf k}}

\newcommand{\R}{{\bf R}}
\newcommand{\Z}{{\bf Z}}
\newcommand{\C}{{\bf C}}
\newcommand{\K}{{\tt K}}
\newcommand{\Naaa}{{\bf N}}
\newcommand{\gaaa}{{\bf g}}
\newcommand{\maaa}{{\bf m}}
\newcommand{\aaaa}{{\bf a}}
\newcommand{\naaa}{{\bf n}}

\newcommand{\cZ}{{\cal Z}}
\newcommand{\cE}{{\cal E}}
\newcommand{\cI}{{\cal I}}
\newcommand{\cC}{{\cal C}}
\newcommand{\mod}{{\rm mod}}
\newcommand{\cK}{{\cal K}}
\newcommand{\cA}{{\cal A}}
\newcommand{\cEp}{{{\cal E}^\prime}}
\newcommand{\cU}{{\cal U}}
\newcommand{\Hom}{{\mbox{\rm Hom}}}
\newcommand{\vol}{{\rm vol}}

\newcommand{\End}{{\mbox{\rm End}}}
\newcommand{\Ext}{{\mbox{\rm Ext}}}
\newcommand{\rk}{{\mbox{\rm rank}}}
\newcommand{\im}{{\mbox{\rm im}}}

\newcommand{\cF}{{\cal F}}
\newcommand{\Ree}{{\rm Re }}

\newcommand{\ee}{{\rm e}}

\newcommand{\tr}{{\mbox{\rm tr}}}

\newcommand{\Ad}{{\mbox{\rm Ad}}}
\newcommand{\ad}{{\mbox{\rm ad}}}

\newcommand{\coker}{{\rm coker}}
\newcommand{\id}{{\mbox{\rm id}}}
\newcommand{\ord}{{\mbox{\rm ord}}}
\newcommand{\nat}{{\bf  N}}
\newcommand{\supp}{{\mbox{\rm supp}}}
\newcommand{\spec}{{\mbox{\rm spec}}}
\newcommand{\Ann}{{\mbox{\rm Ann}}}
\newcommand{\aca}{{\aaaa_\C^\ast}}

\newcommand{\hc}{{{\cal HC}(\gaaa,K)}}

\def\hB{\hspace*{\fill}$\Box$ \newline\noindent}

\newcommand{\Fin}{{\mbox{\rm Fin}}}
\newcommand{\cS}{{\cal S}}
\newcommand{\orig}{{\cal O}}
\def\hB{\hspace*{\fill}$\Box$ \newline\noindent}
\newcommand{\cL}{{\cal L}}

\newtheorem{prop}{Proposition}[section]
\newtheorem{lem}[prop]{Lemma}

\newtheorem{theorem}[prop]{Theorem}
\newtheorem{kor}[prop]{Corollary}

\newtheorem{con}[prop]{Conjecture}

\begin{document}

\maketitle

\tableofcontents
\section{Introduction}
Let $G$ be a connected semisimple Lie group with finite center.
Let $\gaaa$ be the Lie algebra of $G$ and $K$ be its maximal
compact subgroup. By $\hc$ we denote the category
of Harish-Chandra modules of $G$.

The Harish-Chandra module $(\pi,V_{\pi,K})\in \hc$ has natural globalizations
(analytic, smooth, distribution, and hyperfunction vectors)
$$V_{\pi,\omega}\subset V_{\pi,\infty }\subset V_{\pi,-\infty}\subset
V_{\pi,-\omega}\ .$$
Here $V_{\pi,*}$ is a complete locally convex vector space admitting
a continuous $G$-action such that the underlying $(\gaaa,K)$-module
is $V_{\pi,K}$.

Let $\Gamma\subset G$ be a torsion-free
discrete subgroup of finite covolume.
Then a natural problem is to study the cohomology of
$\Gamma$ with coefficients in the different globalizations $V_{\pi,*}$.
The cohomology groups again become locally convex topological vector spaces.
Our work in \cite{bunkeolbrich947} led us to the following conjecture.
\begin{con}\label{grt}
If $\Gamma$ is cocompact, then for all $p\ge 0$ there are isomorphisms
\begin{eqnarray}H^p(\Gamma,V_{\pi,\omega})&=&
H^p(\Gamma,V_{\pi,\infty})\label{yu1}\\
H^p(\Gamma,V_{\pi,-\omega})&=& H^p(\Gamma,V_{\pi,-\infty})\label{entz}\\
H^p(\Gamma,V_{\pi,*})^*&=& H^{\dim(G/K)-p}(\Gamma,V_{\tilde{\pi},-*}),
\:\:*=\omega,\infty\label{yu2}
\end{eqnarray}
and all vector spaces are Hausdorff and finite-dimensional.
If $\Gamma$ has finite covolume, then
$H^\ast(\Gamma,V_{\pi,-\infty})$ is finite-dimensional and Hausdorff.
\end{con}
In \cite{bunkeolbrich947} we settled this conjecture  in the rank-one case
for $*=\omega$ and cocompact $\Gamma$.
In the present paper we prove the conjecture in the case $\rk_\R(G)=1$
for $*=\infty$ admitting non-cocompact $\Gamma$ of finite covolume
(Section \ref{ggkk}).
As indicated in \cite{bunkeolbrich947}, the main result needed
to extend the method of \cite{bunkeolbrich947} to the case $*=\infty$
was the surjectivity of the $B:=\Omega-\mu$ ($\mu\in\C$, $\Omega$ - Casimir of
$G$)
on the space of sections of moderate growth of homogeneous bundles over
$X:=G/K$.
This is proved in Theorem \ref{sur} of the present paper.

We also obtain a more concrete description of the cohomology groups
in terms of automorphic and cusp forms.

Let now $G$ be of general rank and $\Gamma$ be cocompact.
Employing the recent result of Kashiwara-Schmid \cite{kashiwaraschmid94}, Thm.
6.13,
one can show (\ref{yu2}) for $*=\omega$
and $\dim\:H^*(\Gamma,V_{\pi,\pm\omega})<\infty$.
Using Theorem
\ref{casss} below for $*=\infty$ it is possible to prove
(\ref{yu1}) and (\ref{entz}).
We will explain this in a forthcoming paper.
The case of non-cocompact $\Gamma$ remains open for
$\rk_R(G)>1$.

Our main motivation for considering the $\Gamma$-cohomology of
globalizations of Harish-Chandra modules was to prove a conjecture of Patterson
about the singularities of the Selberg zeta functions associated to $\Gamma$.

Let $G$ be a semisimple Lie group of real rank one and
and $P=MAN$ be a minimal parabolic subgroup of $G$.
Let $\aaaa,\naaa$ be the Lie-algebras of $A,N $.
Let $(\sigma,V_\sigma)\in\hat{M}$ be
an irreducible representation of $M$ and $\Gamma\subset G$ be a discrete
cocompact
torsion-free  subgroup.
Then the Selberg zeta function $Z_\Gamma(s,\sigma)$, $s\in\aca$,
is defined as the analytic continuation of the infinite product
$$Z_\Gamma(s,\sigma)=\prod_{[g]\in C\Gamma,n_\Gamma(g)=1}\prod_{k=0}^\infty
\:\det\left(1-a_g^{-\rho-s}S^k(Ad(m_ga_g)_\naaa^{-1})\otimes
\sigma(m_g)\right).$$
Here $C\Gamma$ is the set of hyperbolic conjugacy classes in $\Gamma$,
$n_\Gamma(g)$
is the maximal number $n\in\nat$ such that $g=h^n$ for some $h\in\Gamma$
and $m_g\in M$, $a_g\in A^+$  are such that $g$ is conjugated in $G$ to
$m_ga_g$.
$S^k(\Ad(m_ga_g)_\naaa^{-1})$ stands for the $k$'th symmetric power
of $\Ad(m_ga_g)^{-1}$ restricted to $\naaa$ and $\rho\in\aaaa^\ast$ is defined
by
$\rho(H):=\frac{1}{2}\tr(\ad(H)_\naaa)$.
The infinite product converges for $\Ree(s) > \rho$. In this generality
the Selberg zeta function was introduced by Fried \cite{fried86}.
He applied Ruelle's techniques of hyperbolic dynamics and gave a
meromorphic continuation of $Z_\Gamma(s,\sigma)$.

The parameters $(\sigma,\lambda)\in \hat{M}\times\aca$ also define a principal
series
representation of $G$
on the Hilbert space
$$H^{\sigma,\lambda}=\{f:G\rightarrow
V_\sigma\:|\:f(gman)=a^{\lambda-\rho}\sigma(m)^{-1}f(g),\:\forall man\in MAN,
f_{|K}\in L^2\},$$
where $G$ acts by the left regular representation.
By $H^{\sigma,\lambda}_{-\infty}$ we denote the space of its distribution
vectors.

S. Patterson \cite{patterson93} conjectured  the relationship of the
singularities of Selberg zeta functions
with the  cohomology of $\Gamma$ with coefficients in the distribution vectors
of principal series representations. In \cite{bunkeolbrich947} we have shown an
analog of Pattersons conjecture for the hyperfunction vectors of principal
series representations.
Pattersons original conjecture now follows from (\ref{entz}).
\begin{theorem}\label{umth1}
The cohomology $H^p(\Gamma,H^{\sigma,\lambda}_{-\infty})$ is finite dimensional
for all $p\ge 0$,
\begin{eqnarray}
\chi(\Gamma,H^{\sigma,\lambda}_{-\infty})=\sum_{p=0}^\infty (-1)^p \dim
H^p(\Gamma,H^{\sigma,\lambda}_{-\infty})&=&0,\label{uass1}\\
\chi(\Gamma,\hat{H}^{\sigma,\lambda}_{-\infty})=\sum_{p=0}^\infty (-1)^p \dim
H^p(\Gamma,\hat{H}^{\sigma,\lambda}_{-\infty})&=&0,\label{uass2}
\end{eqnarray}
and the order of $Z_\Gamma(s,\sigma)$ at $s\in\aca$
can be expressed in terms of the group cohomology of $\Gamma$ with coefficients
in $H^{\sigma,\lambda}_{-\infty}$ as follows :
\begin{eqnarray}\ord_{s=\lambda\not=0} Z_\Gamma(s,\sigma)&=&-\sum_{p=0}^\infty
(-1)^p p \dim H^p(\Gamma,H^{\sigma,\lambda}_{-\infty}),\label{pat}\\
\mbox{\em ord}_{s=0} Z_\Gamma(s,\sigma)&=&-\sum_{p=0}^\infty (-1)^p p \dim
H^p(\Gamma,\hat{H}^{\sigma,0}_{-\infty}),
\label{pat0}
\end{eqnarray}
where $\hat{H}^{\sigma,\lambda}_{-\infty}$ is a certain non-trivial extension
of $H^{\sigma,\lambda}_{-\infty}$ with itself.
\end{theorem}
Our description of the cohomology for non-cocompact $\Gamma$ of finite covolume
implies that the analog
of Pattersons conjecture is false in that case.

A very detailed analysis of the cohomology of Fuchsian groups of the first kind
with coefficients in principal series
representations is given in Section \ref{fufu}. We express the cohomology
in terms of automorphic and cusp forms. For the non-irreducible
principal series representations we find explicit formulas
for the dimensions of the cohomology groups in terms of the topology of
$Y=\Gamma\backslash  X$.

The methods of this paper can also be used to study the $\naaa$-cohomology
of globalizations of Harish-Chandra modules.

Let $G$ be of general rank and
$P\subset G$ be a parabolic subgroup
with Langlands decomposition $P=M_PA_PN_P$.
Let $\naaa_P$ be the Lie algebra of $N_P$.

If $(\pi,V_{\pi,K})\in\hc$, then by a theorem of Hecht-Schmid
\cite{hechtschmid83}
the Lie algebra cohomology groups
$H^*(\naaa_P,V_{\pi,K})$ have natural structures of
modules in ${\cal HC}(\maaa_P,K_P)$,
where $K_P=K\cap M_P$ and $\maaa_P$ is the Lie algebra of $M_P$.
Moreover, all elements of  $H^*(\naaa_P,V_{\pi,K})$ are finite
under $\aaaa_P$, where $\aaaa_P$ is the Lie algebra of $A_P$.

It is natural to ask whether globalization is compatible with
$\naaa_P$-cohomology.
The standard $\naaa_P$-cohomology complex induces a natural topology on
$H^p(\naaa_P,V_{\pi,*})$ such that $H^p(\naaa_P,V_{\pi,*})$ becomes a
continuous
representation of  $M_PA_P$.
The natural theorem is
\begin{theorem}\label{casss}
There are topological isomorphisms of $M_PA_P$ modules
\begin{eqnarray*}H^p(\naaa_P,V_{\pi,*})&\cong& H^p(\naaa_P,V_{\pi,K})_*\ , \\
H^p(\naaa_P,V_{\pi,-*})^\ast&\cong&
H^{\dim(\naaa_P)-p}(\naaa_P,V_{\tilde{\pi},K})_{ *}\otimes
\Lambda^{\dim(\naaa_P)}\naaa_P\ ,\quad \forall p\ge 0\end{eqnarray*}
for $*= \infty, \omega$, where $(\tilde{\pi},V_{\tilde{\pi},K})$ is the dual
Harish-Chandra module.
\end{theorem}
For $*=\infty$ this was proved by Casselman (unpublished).
In the rank-one case there were previous partial results by Osborne
\cite{osborne72}.
In the present paper we give a proof in the rank-one case using different
methods (Theorem \ref{haha}).

In \cite{bunkeolbrich947} we were able to prove this theorem in the case
$\rk_\R(G)=1$ for
$*=\omega$.
For $G$ of general rank and $*=\omega$ the theorem
was announced by Bratten and Hecht-Taylor (compare \cite{hechttaylor93}).
It is now an easy consequence of the above-mentioned result
\cite{kashiwaraschmid94}, Thm. 6.13.

\section{Surjectivity on weighted spaces}

{}From now on let $G$ be a connected semisimple Lie
group with finite center of {\em real rank one} and $K$ be its maximal compact
subgroup.
Let $G\times_KV_\gamma=E$ be a $G$-homogeneous bundle over
the rank-one symmetric space $X=G/K$. Let $\cE=C^\infty(X,E)$
be the space of smooth sections of $E$.
We define the  increasing sequence of subspaces $S_R\cE$, $R\in\R$,
by
$$S_R\cE:=\{f\in\cE\:|\: p_{-R,L }(f)<\infty\: \forall L\in\cU(\gaaa) \},$$
where the seminorms $p_{R,L}$ are defined by
$$p_{R,L}(f) :=\sup_G \ee^{R\dist(gK,\orig)}  |f(Lg)|\ , $$
$\orig\in X$ is the class $[K]\in G/K$ considered as the base point of $X$,
and we view $f$ as a function on $G$
with values in $V_\gamma$ satisfying $f(gk)=\gamma(k)^{-1}f(g)$, $\forall k\in
K$.
 $S_R\cE$ is a Fr\'echet space and
for $R^\prime\ge  R$ we have a  continuous
inclusion $S_R\cE\hookrightarrow S_{R^\prime}\cE$.
Let $S_R\cE^\prime:=(S_{-R}\cE)^\ast$ be the topological conjugate dual of
$S_{-R}\cE$.
Then for $R \ge R^\prime$ there is a continuous  projection
$S_{R^\prime}\cE^\prime\rightarrow S_R\cE^\prime$.
We also consider
\begin{eqnarray*}
S_\infty\cE^\prime&:=&\lim_{\stackrel{\longrightarrow}{R}} S_R\cE^\prime\ ,\\
S_\infty\cE&:=&\{f\in\cE\:|\:\forall L\in \cU(\gaaa)\: \exists R\in\R \:s.t.
\:p_{R,L}(f)<\infty\} \ .
\end{eqnarray*}
The space $S_\infty\cE$ is a limit of Fr\'echet spaces in the following way.
Let $\cI$ be the partially ordered set of all monotone functions
$u:\nat_0\rightarrow \R$,
where $u\ge u^\prime$, iff $u(n)\ge u^\prime(n)$, $\forall n\in\nat$.
For any $u\in \cI$ we define the Fr\'echet space
$$S_u\cE:=\{f\in\cE\:|\: p_{-u(\deg(L)),L}(f)<\infty\:\:\:  \forall
L\in\cU(\gaaa) \}\ ,$$
where $\deg(L)$ is the order of $L$ as a differential operator on $G$.
Then
$$S_\infty \cE=\lim_{\stackrel{\longrightarrow}{u\in \cI}} S_u\cE\ .$$
We equip $S_\infty \cE$ with the topology of the direct limit.
The symbol $S_\infty\cE$ is an abuse
of notation since $$S_\infty\cE\not= \lim_{\stackrel{\longrightarrow}{R\in
\R}}S_R\cE\ .$$
$S_\infty\cE$ is called the space of sections of $E$ of moderate growth.

The space $S_\infty\cE^\prime$ is the topological conjugate dual of $\cap_{R}
S_R\cE$.
Since $C_c^\infty(X,E)$ is dense in $\cap_{R} S_R\cE$,
we have an inclusion $S_\infty\cEp\hookrightarrow C^{-\infty}(X,E)$.

Let $\Omega$ be the Casimir operator of $G$ and for $\mu\in\C$ let
$B:=\Omega-\mu$.
Let $\C[B]$ be the ring of all polynomials in $B$.
We consider the functor $\Fin_\mu$ on the category of $\C[B]$-modules
defined by
$$\Fin_\mu(V)=\{v\in V\:|\: \exists l\ge 0 \quad B^lv=0  \}\ .$$
This functor is left exact and its higher derived functors are denoted by
$\Fin_\mu^i(V)$.
Since the inductive limit  is an exact functor on the category of
$\C[B]$-modules
and $$\Fin_\mu(V)\cong\lim_{\stackrel{\longrightarrow}{j}} \Hom_{\C[B]}(
\C[B]/(B )^j,V)$$
we have
\begin{equation}\label{iop}\Fin_\mu^i(V)\cong
\lim_{\stackrel{\longrightarrow}{j}} \Ext^i_{\C[B]}(\C[B]/(B )^j,V)\
.\end{equation}
Since $\C[B]$ is a regular ring of dimension one
we have $\Fin_\mu^i(V)=0$ for  $i\ge 2$.
A $\C[B]$-module  $V$ is called $\Fin_\mu$-acyclic, iff $\Fin_\mu^1(V)=0$.
\begin{theorem}\label{fin}
$S_\infty\cEp$ is $\Fin_\mu$-acyclic.
\end{theorem}
\proof
We consider the Schwartz space
$$S_{-log}\cE:=\{f\in\cE\:|\: q_{L,k}(f)<\infty \quad \forall k\in\nat , L\in
\cU(\gaaa)\}\ ,$$
where $$q_{L,k}(f)^2=\int_G \dist(gK,\orig)^k |f(Lg )|^2 dg$$
and its topological conjugate dual space
$$S_{log}\cEp=(S_{-log}\cE)^\ast\ .$$
\begin{lem}\label{ssur}
The inclusion $S_{log}\cEp \hookrightarrow S_\infty\cEp$
induces a surjection
$$\Fin^1_\mu(S_{log}\cEp)\rightarrow \Fin_\mu^1(S_\infty\cEp)\ .$$
\end{lem}
\proof
Let $V$ be  a $\C[B]$-module.
We can compute
$\Ext^1_{\C[B]}(\C[B]/(B )^k,V)$ using the Koszul complex
$$0\rightarrow \Ext^0_{\C[B]}(\C[B]/(B )^k,V)\rightarrow
V\stackrel{B^k}{\rightarrow}V\rightarrow
\Ext^1_{\C[B]}(\C[B]/(B )^k,V)\rightarrow 0\ .$$
In order to perform the direct limit in (\ref{iop})
note that the map
$$\Ext^1_{\C[B]}(\C[B]/(B^k),V)\rightarrow
\Ext^1_{\C[B]}(\C[B]/(B^{k+1}),V)$$
is induced by $B:V\rightarrow V$.

Hence Lemma \ref{ssur} follows from the next result.
\begin{lem}\label{iter}
Let $k\ge 1$.
For any $f\in S_\infty\cEp$ there  exists a $g\in S_\infty\cEp$  with
$f-B^k g\in S_{log}\cEp$.
\end{lem}
In fact, Lemma \ref{iter} implies
that any element $f\in S_\infty\cEp$ representing some class in\linebreak[4]
$\Ext^1_{\C[B]}(\C[B]/(B )^k,S_\infty\cEp)$  can be replaced by $f-B^kg\in
S_{log}\cEp$
representing the same class in $\Ext^1_{\C[B]}(\C[B]/(B )^k,S_\infty\cEp)$.\\
\proof
It is enough to show that for $f\in S_R\cEp$ there exists
$g\in S_R\cEp$ such that $f-B^kg\in S_{R-1}\cEp$.
In fact, it is sufficient that this holds for all weights $R$
in a dense subset of $\R$.
Then the assertion of the lemma follows by a finite iteration.

We consider polar coordinates $ \aaaa^+\times K/M $ of $X\setminus \{\orig\}$.
We identify the fibre of $E$ over $(a,kM)$ with the fibre over $(a^\prime,kM)$
using the radial parallel transport.

There is a constant coefficient
operator $B^k_{rad}$ on $\aaaa^+$ such that  $B^k-B^k_{rad}=\ee^{-a}Q$ (we
identify
$\aaaa^+$ with $\R_+$ such that the short root has length one) and
$Q$ has bounded coefficients up to infinity.
Note that $Q$ also contains differentiations along $\aaaa^+$ but with
exponentially decreasing
prefactors.
There are an endomorphism $F(kM)\in\End(E_{(a,kM)})$ and $c\in \C$ such that
$$B_{rad} = - \frac{d^2}{da^2}+c\frac{d}{da}+F(kM)\ .$$
We decompose the bundle $E_{|X\setminus \orig}=\sum_{\sigma}E(\sigma)$
according to the
eigenvalues of $F(kM)$ (which are independent of $kM$).
There are  $x_\sigma,y_\sigma\in \C$ such that
$B_{rad}=\sum_{\sigma}B_{rad}(\sigma)$ and
$B_{rad}(\sigma)=-(\frac{d}{da}-\bar{x_\sigma})(\frac{d}{da}-\bar{y_\sigma}) $.
We exclude the weights $R= \Ree(x_\sigma), \Ree(y_\sigma)$ for all $\sigma$
from the following consideration.

Let $S_R\cE_0\subset S_R\cE$ be the subspace of all sections vanishing in a
neighbourhood
of the origin $\orig\in X$. Then $S_R\cE_0$ is a limit of Fr\'echet spaces. Let
$S_R\cEp_0$
be the topological conjugate dual of $S_R\cE_0$. Then $S_R\cEp\subset
S_R\cEp_0$.
Choose a cut-off function $\chi\in C^\infty(\aaaa^+)$ such that $\chi(a)=0$ for
$a<1$ and $\chi(a)=1$ for $a>2$. Then the multiplication by $\chi$ defines a
continuous
operator $S_R\cEp_0 \rightarrow S_R\cEp$.

Below we will construct a continuous solution operator
$H:S_R\cE_0\rightarrow S_R\cE_0$ such that
$ H(B^*_{rad})^k -\id : S_R\cE_0\rightarrow  C_c^\infty(X\setminus\orig,E)$
is continuous and $\supp (H(B^*_{rad})^k -\id)f \subset [1,2]\times K/M$,
$\forall f\in S_R\cE_0$.
Here $(B_{rad}^*)^k $ is the formal adjoint of $B^k_{rad}$.
The adjoint $H^*:S_R\cEp_0\rightarrow S_R\cEp_0$ then has the property that
$B^k_{rad}H^*\phi-\phi\in  \cEp$ for $\phi\in S_R\cEp_0$.
Here $\cEp$ is the space of distributions with compact support.
The composition $\chi\circ H^*_{|S_R\cEp}:S_R\cEp\rightarrow S_R\cEp$
satisfies $B^k_{rad}\chi\circ H^*\phi-\phi\in \cEp$, $\phi\in S_R\cEp$. Hence
$B^k H^*\phi-\phi=e^{-a}QH^*\phi\:(\mod\: \cE^\prime)$.
But $e^{-a}QH^*\phi\in S_{R-1}\cEp$.
This proves the lemma assuming that we have already constructed $H$.

Note that $ B^*_{rad}(\sigma)^k$
is a product of operators of first order
 $\frac{d}{da}+x_\sigma$, $\frac{d}{da}+y_\sigma$.
It  is enough to construct  solution operators
$H_{x_\sigma},H_{y_\sigma}:S_R\cE_0(\sigma)\rightarrow S_R\cE_0(\sigma)$
such that $\supp (H_{x_\sigma}(\frac{d}{da}+x_\sigma)f-f)\subset [1,2]\times
K/M$, $\forall f\in S_R\cE_0(\sigma)$,
and similarly for $H_{y_\sigma}$.
Then $H:=(-1)^k\sum_{\sigma}H_{x_\sigma}^kH_{y_\sigma}^k$ has the required
properties.

 For $R-\Ree(x_\sigma)>0$ we set (simplifying the notation by omitting the
angular variable)
$$(H_{x_\sigma}f)(a):=-\chi(a)\ee^{-x_\sigma a}\int_a^\infty \ee^{x_\sigma b}
f(b) db\ .$$
Since $f\in S_R\cE_0(\sigma)$  we have for some $C<\infty$ and all
$b\in\aaaa^+$  that
$|f(b)|\le C \ee^{-Rb}$. Hence the integral converges.
If $R-\Ree(x_\sigma)< 0$, then we set
$$(H_{x_\sigma}f)(a):=\chi(a)\ee^{-x_\sigma a}\int_0^a \ee^{x_\sigma b}   f(b)
db\ .$$
It is easy to check that $\supp (H_{x_\sigma} (\frac{d}{da}+x_\sigma)
f-f)\subset [1,2]\times K/M$.
We must show that $H_{x_\sigma}:S_R\cE_0(\sigma)\rightarrow S_R\cE_0(\sigma)$
and that $H_{x_\sigma}$ is continuous.

To define the topology of $S_R\cE_0(\sigma)$ it is sufficient to consider the
set of seminorms
$\tilde{p}_{R,L}$ with $L=(L_1,L_2)\in\cU(\aaaa)\times \cU(\kaaa)$.
Set $(Lf)(a,kM)=f(L_1a,L_2kM)$.
Then $$\tilde{p}_{R,L}: =\sup_{a\in\aaaa^+, k\in K} \ee^{Ra} |(Lf)(a,kM)| \ .$$
It  is clear that $LH_xf=H_x Lf + W_{L}f$, where
$W_{L}:S_R\cE_0(\sigma)\rightarrow C_c^\infty(X\setminus \orig,E(\sigma))$ is
continuous.
Thus in order to show that $H_{x_\sigma}$ is continuous it is enough
to verify the estimate
$$\tilde{p}_{R,1}(H_{x_\sigma}f) \le C \tilde{p}_{R,1}(f)\ .$$
We employ that $|f(a,kM)|\le \tilde{p}_{R,1}(f)\ee^{-Ra}$.
If $R-\Ree(x_\sigma)>0$, then
\begin{eqnarray*} \tilde{p}_{R,1}(H_{x_\sigma}f) &=&
\sup_{a\in\aaaa^+,k\in K} \ee^{Ra}|(H_{x_\sigma}f)(a,kM)|\\  &\le &
\sup_{a\in\aaaa^+,k\in K} \ee^{Ra}\tilde{p}_{R,1}(f) \ee^{-\Ree(x_\sigma)a}
\int_{a}^\infty \ee^{(\Ree(x_\sigma)-R)b}db\\
&\le & C \tilde{p}_{R,1}(f)\ .
\end{eqnarray*}
If $R-\Ree(x_\sigma)<0$, then
\begin{eqnarray*}\tilde{p}_{R,1}(H_{x_\sigma}f) &=&
\sup_{a\in\aaaa^+,k\in K} \ee^{Ra}|(H_{x_\sigma}f)(a,kM)|\\  &\le &
\sup_{a\in\aaaa^+,k\in K} \ee^{Ra}\tilde{p}_{R,1}(f)\ee^{-\Ree(x_\sigma)a}
\int_{1}^a \ee^{ (\Ree(x_\sigma)-R)b}db\\
&\le & C \tilde{p}_{R,1}(f)\ .
\end{eqnarray*}
This finishes the construction of $H$ and hence  the proof of  Lemma
 \ref{ssur}.\hB
\begin{lem}\label{logfin}
$S_{log}\cEp$ is $\Fin_\mu$-acyclic.
\end{lem}
\proof
Using the generalization of the Helgason-Fourier transform to bundles
(see \cite{bunkeolbrich955}, 1.4.3 , Branson-Olafsson-Schlichtkrull
\cite{bransonolafssonschlichtkrull92})
and a result of Arthur \cite{arthur75}
we can identify the Schwartz space
$S_{-log}\cE$ with
$$[\cS(\imath \aaaa^\ast)\hat{\otimes} C^\infty(K\times_M V_\gamma)]^W\oplus
\oplus_{i=1}^r V_{i,\infty } \ .$$
Here $V_{i,\infty}$ are the spaces of smooth vectors   of certain irreducible
discrete series representations, $\cS(\imath\aaaa^\ast)$
is the Schwartz space on $\imath\aaaa^\ast$ and $W=\Z_2$ is the Weyl group.
$\hat{\otimes}$ denotes the nuclear tensor product.
The operation of $W$ on $\cS(\imath \aaaa^*)\hat{\otimes} C^\infty(K\times_M
V_\gamma)$
is implemented by a family of unitary operators
$$A_\lambda:C^\infty(K\times_M V_\gamma)\rightarrow C^\infty(K\times_M
V_\gamma),\quad \lambda\in\imath \aaaa^*\ ,$$
which are closely related to
 Knapp-Stein intertwining operators of principal series representations.
The non-trivial element $w\in W$ acts by
$(wf)(-\lambda,k)=(A_\lambda f(\lambda,.))(k)$.
 Dually, we can identify
$S_{log}\cEp$ with
$$[\cS^\prime(\imath \aaaa^*)\hat{\otimes} C^{-\infty}(K\times_M
V_\gamma)]^W\oplus \oplus_{i=1}^r V_{i,-\infty}^\prime \ .$$

On the Fourier image the operator $B$ acts as multiplication
by the polynomial $P(\lambda):=\lambda^2+\rho^2-\gamma(\Omega_M)-\mu$,
$\lambda\in\imath\aaaa^*$, on the first
summand and as a scalar on $V_{i,-\infty}^\prime$.

It follows that $\Fin_\mu^1(V_{i,-\infty}^\prime)=0$.
In order to see that
$$\Fin_\mu^1([\cS^\prime(\imath \aaaa^*)\hat{\otimes} C^{-\infty}(K\times_M
V_\gamma)]^W)=0$$
we show that the multiplication  by $P$
is surjective on   $[\cS^\prime(\imath \aaaa^*)\hat{\otimes}
C^{-\infty}(K\times_M V_\gamma)]^W$.
We  consider the dual situation and show that the multiplication by
$P$ on $[\cS(\imath \aaaa^\ast)\hat{\otimes} C^\infty(K\times_M V_\gamma)]^W$
is injective and has a closed range.
Injectivity is easy since $P(\lambda)$ is non-invertible at most on a set of
codimension $1$
in $\imath \aaaa^*\times K/M$.
Note that $C^\infty(K\times_M V_\gamma)=\oplus_{\sigma\in \hat{M}}
\oplus_{h=1}^{[\gamma:\sigma]}
C^\infty(K\times_M V_{\sigma})$
and $P(\lambda)$ acts by the scalar polynomial
$P_\sigma(\lambda)=\lambda^2+\rho^2-\sigma(\Omega_M)-\mu$
on each summand $\cS^\prime(\imath \aaaa^*)\hat{\otimes} C^\infty(K\times_M
V_{\sigma})$.
Let $Z_\sigma$ be the set of zeros of $P_\sigma(\lambda)$.
Let $n_\sigma(\lambda)$ be the order of the zero of $P_\sigma$ at $\lambda \in
Z_\sigma$.
Then the range of the multiplication by $P_\sigma$ on $\cS(\imath
\aaaa^\ast)\hat{\otimes} C^\infty(K\times_M V_{\sigma})$ consists of all
sections
$f\in  \cS(\imath \aaaa^\ast)\hat{\otimes} C^\infty(K\times_M V_{\sigma}) $
such that $f(\lambda,kM)$ has a zero of order $n_\sigma$ for all $\lambda \in
Z_\sigma$  and for all $k\in K$.
This vanishing condition is a closed condition.
The range of $P$ on $[\cS(\imath \aaaa^\ast)\hat{\otimes} C^\infty(K\times_M
V_\gamma)]^W$
consists of all $f\in [\cS(\imath \aaaa^\ast)\hat{\otimes} C^\infty(K\times_M
V_\gamma)]^W$
such that $f_\sigma(\lambda,kM)$ has a zero of order $n_\sigma(\lambda)$ for
all $\lambda\in Z_\sigma$, $k\in K$ and $\sigma$.
This condition is again closed.

This finishes the proof of Lemma \ref{logfin} and hence of Theorem
\ref{fin}.\hB
 \begin{theorem}\label{sur}
The operator $B:S_\infty\cE\rightarrow S_\infty\cE$
is surjective.
\end{theorem}
\proof
Let $S_\infty\cEp(B)$ and $C$ be defined by the following exact sequence:
\begin{equation}\label{mmnnj}0\rightarrow S_\infty\cEp(B)\rightarrow
S_\infty\cEp\stackrel{B}{\rightarrow} S_\infty \cEp\rightarrow C\rightarrow 0\
.\end{equation}
\begin{lem}\label{dissur}
$C=0$.
\end{lem}
\proof
We have  $BC=0$ and $BS_\infty\cEp(B)=0$. Hence
$\Fin_\mu^1(S_\infty\cEp(B) )=0$, $\Fin_\mu(S_\infty\cEp(B) )=S_\infty\cEp(B)
$, $\Fin_\mu^1(C)=0$, and
$\Fin_\mu(C)=C$.
Thus applying $\Fin_\mu$ to the exact sequence of $\Fin_\mu$-acyclic spaces
 (\ref{mmnnj}) we obtain the exact sequence
$$ 0\rightarrow S_\infty\cEp(B)  \rightarrow \Fin_\mu
(S_\infty\cEp)\stackrel{B}{\rightarrow} \Fin_\mu(S_\infty \cEp)\rightarrow
C\rightarrow 0\ .$$
Let $S_\infty\cEp(B^k) :=\ker(B^k:S_\infty\cEp\rightarrow S_\infty\cEp)$. Then
$$\Fin_\mu(S_\infty\cEp)=\lim_{\stackrel{\longrightarrow}{k}}S_\infty\cEp(B^k)
$$
and in order to show that $C=0$ it is enough to  show that
\begin{equation}\label{loo}
\lim_{\stackrel{\longrightarrow}{k}}\coker(B:S_\infty\cEp(B^k)\rightarrow
S_\infty\cEp(B^k))=0\ .
\end{equation}
 The space of $K$-finite vectors $S_\infty\cEp(B^k)_K$
is an admissible $(\gaaa,K)$-module and $S_\infty\cEp(B^k)$
is its canonical distribution vector globalization (\cite{wallach92}, Ch. 11).
Since the globalization functor is exact (Wallach \cite{wallach92}, Ch. 11,
Casselman \cite{casselman89})
to show (\ref{loo}) is suffices to verify that
\begin{equation}\label{huhu}
\lim_{\stackrel{\longrightarrow}{k}}\coker(B:S_\infty\cEp(B^k)_K\rightarrow
S_\infty\cEp(B^k)_K)=0 \ .
\end{equation}
Since the multiplication by $B$ on
$(\cU(\gaaa)\otimes_{\cU(\kaaa)}V_{\tilde{\gamma}})_K$ is injective we obtain
by $K$-type wise algebraic dualization
$$\coker(B:\Hom_{\cU(\kaaa)}(\cU(\gaaa),V_\gamma)_K\rightarrow
\Hom_{\cU(\kaaa)}(\cU(\gaaa),V_\gamma)_K)=0\ .$$
Thus $\Hom_{\cU(\kaaa)}(\cU(\gaaa),V_\gamma)_K$ is $\Fin_\mu$-acyclic.
Applying $\Fin_\mu$ to the exact sequence of $\Fin_\mu$-acyclic spaces
$$
0\rightarrow \Hom_{\cU(\kaaa)}(\cU(\gaaa),V_\gamma)_K(B)\rightarrow
\Hom_{\cU(\kaaa)}(\cU(\gaaa),V_\gamma)_K
\stackrel{B}{\rightarrow}\Hom_{\cU(\kaaa)}(\cU(\gaaa),V_\gamma)_K \rightarrow 0
$$
we obtain the exact sequence
\begin{eqnarray*}
\lefteqn{0\rightarrow \Hom_{\cU(\kaaa)}(\cU(\gaaa),V_\gamma)_K(B)\rightarrow
\Fin_\mu( \Hom_{\cU(\kaaa)}(\cU(\gaaa),V_\gamma)_K)
\stackrel{B}{\rightarrow}}\\
&&\hspace{5cm}\Fin_\mu(\Hom_{\cU(\kaaa)}(\cU(\gaaa),V_\gamma)_K) \rightarrow 0\
{}.
\end{eqnarray*}
The exactness at the last place  is  equivalent to
$$
\lim_{\stackrel{\longrightarrow}{k}}\coker(B:\Hom_{\cU(\kaaa)}
(\cU(\gaaa),V_\gamma)_K(B^k)\rightarrow
\Hom_{\cU(\kaaa)}(\cU(\gaaa),V_\gamma)_K(B^k))=0 \ .
$$
The claim (\ref{huhu}) now follows from
$S_\infty\cEp(B^k)_K\cong \Hom_{\cU(\kaaa)}(\cU(\gaaa),V_\gamma)_K(B^k)$
(compare \cite{bunkeolbrich947}, Lemma 3.5).
\hB
We have shown that $B:S_\infty \cEp\rightarrow S_\infty \cEp$
is surjective.
The $G$-invariant Hermitian scalar product on $E$ induces an inclusion
$S_\infty\cE\subset S_\infty\cEp$.
Theorem \ref{sur} is now a consequence
of the following regularity result.
\begin{lem}\label{reg}
If $f\in S_\infty\cEp$, $Bf\in S_\infty\cE$, then $f\in S_\infty\cE$.
\end{lem}
\proof
Let $\tilde{W}$ be a $G$-invariant fundamental solution of $B$, i.e.,
$B\tilde{W}=\tilde{W}B=\id$ on $C_c^\infty(X,E)$.
Consider $X\times X$ with the diagonal $G$-action.
Let $\chi\in {}^GC^\infty(X\times X)$ be a $G$-invariant cut-off
function defined by $\chi(x,y)=\rho(\dist(x,y))$, where $\rho\in C^\infty(\R)$,
$\rho(r)\in [0,1]$,
$\rho(r)=1$ for $r<1$ and $\rho=0$ for $r>2$.

Multiplying the distributional kernel of $\tilde{W}$ by $\chi$
we obtain a $G$-invariant parametrix $W$ of $B$.
It has finite propagation
(the support of $W\phi$ is contained in a $1$-neighbourhood
of the support of $\phi$ for all $\phi\in C_c^\infty(X,E)$)
and is applicable to arbitrary distributions in $C^{-\infty}(X,E)$.
Let $f \in S_\infty\cEp$ and $Bf=F\in S_\infty\cE$.
Then $f=WF+Sf$, where $S$ is a $G$-invariant smoothing
operator on $E$.
We must show that $WF,Sf\in S_\infty\cE$.

Let $L\in \cU(\gaaa)$. Then
$$(Sf)(Lg) = \langle  f ,s(LgK,. ) \rangle \ ,$$
where $s$ is the integral kernel of $S$ and
$s(LgK,.)\in C_c^\infty(X, E)$.
For any $L_1\in \cU(\gaaa)$ and $R\in \R$ there is a $P\in  \R$ depending on
$L_1,L$ such that
 $$p_{R,L_1}(s(LgK, .)) = \sup_{h\in G} \ee^{R\dist(\orig,hK)} |s(LgK,L_1hK)|
\le C \ee^{(P+R)\dist(\orig,gK)}\ ,$$
where $C$ may depend on $L ,L_1,R,P$ but not on $g$.
In fact, for all $A,B\in \cU(\gaaa)$ we have $\sup_{g,h\in G}|s(gA,hB)|<
\infty$ since
$S$ is $G$-invariant.
If $f\in S_R\cEp$ for some $R$, then there is a finite set
$\{L_1,\dots,L_r\}\subset \cU(\gaaa)$
such that for all $L\in \cU(\gaaa)$ there exists
$P\in \R$ with
\begin{eqnarray*}
p_{-(R+P),L}(Sf) &=& \sup_{g\in G} \ee^{-(R+P)\dist(\orig,gK)} \langle f,
s(LgK,. ) \rangle \\
&\le& C(f) \sup_{g\in G, i=1,\dots, r}
\ee^{-(R+P)\dist(\orig,gK)}p_{R,L_i}(s(LgK, .))\\&<&\infty \ .\end{eqnarray*}
This proves $Sf \in S_\infty\cE$.

We now show that $WF\in S_\infty\cE$.
Let $L\in \cU(\gaaa)$. We must show that for some $R\in\R$
$$p_{R,L}(WF)=\sup_{g\in G} \ee^{R\dist(\orig,gK)} |(WF)(LgK)|<\infty\ .$$
Let $w(gK,hK)$ be the distributional kernel of $W$.
The family $gK\mapsto w(LgK,.)$ is a smooth family of distributions
such that  $\supp(w(LgK,.) )$ is contained in the unit ball in $X$ with center
at $gK$.
Thus we can write
$$(WF)(LgK)=\langle w(LgK,.) , F\rangle\ .$$
Since $W$ is $G$-invariant there is a finite set
$\{L_1,\dots,L_r\}\subset \cU(\gaaa)$ such that we have
$$|(WF)(LgK)|\le \ee^{P\dist(\orig,gK)} \sup_{h\in G, \dist(hK,gK)\le 1,
i=1,\dots, r}| F(hL_iK)| $$
for some $P\in R$ depending on $L$.
Hence there is another finite set
$\{L^\prime_1,\dots,L^\prime_{r^\prime}\}\subset \cU(\gaaa)$  and an exponent
$Q\in \R$ such that
\begin{eqnarray*}|(WF)(LgK)|&\le& \ee^{Q\dist(\orig,gK)} \sup_{h\in G,
\dist(hK,gK)\le 1,  i=1,\dots, r^\prime}|F(L_i^\prime hK)|\\
&\le&
\ee^{(Q-R)\dist(\orig,gK)}\sup_{i=1,\dots,
r^\prime}p_{R,L_i^\prime}(F)\end{eqnarray*}
Hence $p_{-(Q-R),L}(WF)<\infty$ if $R$ is small enough.

We conclude that $WF\in S_\infty\cE$.
This finishes the proof of the Lemma and of Theorem \ref{sur}.
\hB

  \section{De Rham complexes}\label{sec15}
In this section we prove the local acyclicity
of the weighted de Rham complex on $X$ twisted with
the functions of moderate growth on $G$,
the global acyclicity of the weighted de Rham complex on $X$ twisted with the
functions of moderate growth on $\Gamma\backslash G$, and
the $\naaa$-acyclicity of $S_\infty\cE$ for a homogeneous vector bundle
$E\rightarrow X$.
A complex is called acyclic if it is exact in all non-zero degrees.
\begin{lem}\label{amdt}
For  any vector bundle $E\rightarrow X$
the space $S_\infty\cE$ is $\naaa$-acyclic.
\end{lem}
\proof
Fix a basis $\{X_i\}_{i=1}^{\dim(\naaa)}$ of $\naaa$.  Let $\{X^i\}$ be the
basis of $\naaa^*$ dual to $\{X_i\}$.
Let $Xn\in T_nN$, $X\in\naaa$, $n\in N$, be the fundamental vector at $n$
corresponding to
$-X$.

Let $I_p$ be the set of all multi-indices $\{1\le i_1<i_2<\dots < i_p\le
\dim(\naaa)\}$.
For $I\in I_p$ we define $X^In\in \Lambda^pT^*_nN$ by
$X^In(X_Jn)=\delta_{IJ}$, $\forall J\in I_p$. Here $X_Jn=X_{j_1}n\wedge\dots
\wedge X_{j_p}n$.
Now identify $N\times A\stackrel{\sim}{\rightarrow}X$, $(n,a)\mapsto naK$.
We also employ the $\naaa$-equivariant trivialization $G\times_K
V_\gamma=E\stackrel{\sim}{\rightarrow}N\times A\times V_\gamma$,
$[nak,v]\mapsto (n,a)\times \gamma^{-1}(k)v$.
We identify $C^\infty(N,\Lambda^*T^\ast N\otimes  C^\infty(A)\otimes V_\gamma)$
with $\Lambda^*\naaa^*\otimes \cE$  such that
$\omega\in  C^\infty(N,\Lambda^pT^\ast N\otimes  C^\infty(A)\otimes V_\gamma)$
corresponds
to $(X_I,na)\to \omega(na)(X_In)$, $I\in I_p$, $n\in N$, $a\in A$,
$\omega(na)(X_In)\in V_\gamma$.
Under this identification the $\naaa$-cohomology
complex of $\cE$ becomes the de Rham complex of $N$
twisted with $C^\infty(A)\otimes V_\gamma$.

Using the trivialization $\{X^In\}_{I \in I_p}$  of the $p$-form
bundle we write \linebreak[4]$\omega(na)=\sum_{I\in I_p} \omega_I(na) X^In$.
 The  subspace $\Lambda^p\naaa\otimes S_\infty \cE$ is identified with
\begin{eqnarray*}
 \Lambda^p\naaa\otimes S_\infty \cE&=&\{\omega\in C^\infty(N,\Lambda^pT^\ast
N\otimes  C^\infty(A)\otimes V_\gamma)\\&&|\:\: \forall  L_1\in\cU(\naaa),
L_2\in \cU(\aaaa)\:\:\: \exists R\in \R\\
&& s.t. \sup_{na\in NA, I\in I_p} \ee^{-R\dist(\orig,naK)}
|\omega_I(L_1n,L_2a)| < \infty\}\ .
\end{eqnarray*}
We  define the contraction $\Psi_t$, $t\in[0,1]$, of $N$ by
$\Psi_t(n):=\exp(t\log(n))$.
Let $T_t:=\frac{d}{ds}_{|s=0}\Psi_{t+s}(n)\in T_{\Psi_t(n)}N$.
We set $H_t\omega:=\Psi_t^*(i_{T_t}\omega)$
and $H=\int_0^1 H_t dt$.
Then $$H:C^\infty(N,\Lambda^pT^\ast N\otimes  C^\infty(A)\otimes V_\gamma)
\rightarrow C^\infty(N,\Lambda^{p-1}T^\ast N\otimes  C^\infty(A)\otimes
V_\gamma)$$
 is a zero homotopy of the de Rahm complex
$(C^\infty(N,\Lambda^*T^\ast N\otimes  C^\infty(A)\otimes V_\gamma),d)$.
In order to prove the Lemma we have to show
that this zero homotopy is compatible with the subspaces
$ \Lambda^*\otimes S_\infty \cE$. It is enough to show that
$$t\to H_t\in \Hom^{cont}_\C(\Lambda^p\otimes S_\infty \cE,\Lambda^{p-1}\otimes
S_\infty \cE)$$ is continuous.
Here we equip $\Hom^{cont}$ with the strong topology
(pointwise convergence).
We call a function $f$ on $N$ a polynomial, if $f(\exp(n))$ is a polynomial
on $\naaa$. Using the fact that $N$ is nilpotent, one can easily show that
$$(H_t\omega)(n)=\sum_{J\in I_{p-1},I\in
I_p}\Phi_{I,J}(t,n)\omega_I(\Psi_t(n))X^Jn\ , $$
where $n\to \Phi_{I,J}(t,n)$ are polynomial functions on $N$ (and in fact
also polynomials in $t$).
Let $L_r$ be the set of
tuples $\{i_1,\dots,i_{\dim(\naaa)}\}$, $i_k\in \nat_0$, with
$\sum_{k=1}^{\dim(\naaa)}i_k\le r$.
Let $X(l):=X_1^{l_1}\dots X_{\dim(\naaa)}^{l_{\dim(\naaa)}}\in \cU(\naaa)$,
$l\in L_r$.
By $X(l)n$ we denote the corresponding differential operator $(X_1n)^{l_1}
\dots (X_{\dim(\naaa)}n)^{l_{\dim(\naaa)}}$.
If $L_1\in \cU(\naaa)$, $\deg(L_1)=r$, $L_2\in \cU(\aaaa)$, then also
$$(H_t\omega)_J(L_1n,L_2a)=\sum_{l\in L_r,I\in I_p} \Phi_{l,I,J}(t,n)
\omega_I(X(l)\Psi_t(n),L_2a)$$
with polynomial functions $\Phi_{l,I,J}(t,n)$.
We obtain  the estimate
$$\sup_{J\in I_{p-1}}|(H_t\omega)_J(L_1n,L_2a)|\le C (1+|\log(n)|)^P \sup_{l\in
L_r,I\in I_p} |\omega_I(X(l)\Psi_t(n),L_2a)|\ ,$$
where $P$ is  sufficiently large.
Note that $\ee^{Q\dist(\orig,naK)} \ge c (1+|\log(n)|)^P$ for some $Q,c>0$
and all $na\in NA$.
We also have $\dist(\orig,naK)\ge \dist(\orig,\Psi_t(n)aK)$ for all $t\in
[0,1]$.
Fix $R\in \R$.
Then for all $\omega$ with
$$\sup_{l\in L_r,na\in NA, I\in I_p} \ee^{-R\dist(\orig,naK)}
|\omega_I(X(l)n,L_2a)|=:S(\omega) < \infty$$
we have $$ \sup_{na\in NA, J\in I_{p-1}} \ee^{-(R+Q)\dist(\orig,naK)}
|(H_t\omega)_J(L_1n,L_2a)| < CS(\omega)$$
with $C<\infty$ independent of $\omega$.
These estimates imply that if $\omega\in \Lambda^p\naaa\otimes  S_u\cE$, $u\in
\cI$,
then $H_t\omega \in  \Lambda^{p-1}\naaa\otimes  S_{u+Q}\cE$ and depends
continuously on $t\in[0,1]$.
This proves the required continuity of $H_t$.
Thus $H$ provides a zero homotopy of the $\naaa$-cohomology
complex for $S_\infty \cE$, and the lemma is proved.
 \hB

Fix an Iwasawa decomposition $G=NAK$ and identify  $X\cong NA$. Using these
coordinates
we attach a boundary $\partial X=N\times\{\infty\}$
to $X$ obtaining $\bar{X}$. Note that this boundary is different from
the geodesic boundary.
The motivation for this definition is that small  neighbourhoods
of $x\in\partial X$ look like small neighbourhoods
of  points in the boundary of the Borel-Serre compactification
of quotients $Y=\Gamma\backslash X$ of finite volume.

We want to twist the de Rham complex of $X$ with the space of functions
of moderate growth on $G$. Let $\pi:X\times G\rightarrow X$ denote the
projection onto
the first factor.
Let $L^p:=\Lambda^pT^*X$.
We define
\begin{eqnarray*}
S_\infty\cL^p[G]&:=&\{\omega\in C^\infty(X\times G,\pi^*L^*)\\
&&|\:\forall L_1,L_2\in \cU(\gaaa)\:\:\: \exists R,Q\in\R\\
&&s.t.\:\: \sup_{g,h\in G}\ee^{-R\dist(\orig,gK)}\ee^{-Q\dist(\orig,hK)}
|\omega(L_1g,L_2h)|<\infty
\} \ ,
\end{eqnarray*}
where we consider $\omega$ as a function from $G\times G$ with values in
$\Lambda^pT_\orig^*X $.

There is a natural sheafification $\underline{S_\infty\cL^\ast[G]}$ on
$\bar{X}$
such that for an open set $U\subset \bar{X}$ the vector space
$\underline{S_\infty\cL^\ast[G]}(U)$ consists of those forms $\omega\in
C^\infty(U\times  G,\pi^* L^\ast)$
which satisfy
$$\sup_{g\in U,h\in G}
|\omega(L_1g,L_2h)|\ee^{-R\dist(\orig,gK)}\ee^{-Q\dist(\orig,hK)} <\infty$$
for any $L_1,L_2\in \cU(\gaaa)$ and appropriate $Q,R\in \R$ which may depend on
$L_1,L_2$.
Let $d$ be the differential of the de Rham complex acting trivially with
respect to the second variable $g\in G$.
\begin{lem}\label{mmarrt}
The complex of sheaves $(\underline{S_\infty\cL^\ast[G]},d)$
is acyclic.
\end{lem}
\proof
Let $x\in \bar{X}$. Then we have to show that the complex of germs
$(\underline{S_\infty\cL^\ast[G]}_x,d)$  is acyclic.
If $x\in X$, then we employ the standard homotopy formula
associated to the radial contraction of small balls in $X$ with center at $x$.
We leave to the reader to verify that the zero homotopy is compatible with the
growth conditions with respect to second variable.
We discuss a similar problem in the proof of Lemma \ref{mglobal}.

It remains to consider $x\in \partial X$. Without loss of generality
we can  assume $x=e\times\{\infty\}$, where $e\in N$ is the identity.
Let $U_i\subset \naaa$, $i\in \nat$, be a fundamental
system  of balls around $0$. Then $W_i:=\exp(U_i) (i,\infty]\subset NA$
is a fundamental system of neighbourhoods of $x$.
Here we identified $\aaaa\cong\R$ such that
$(i,\infty]=\exp((i,\infty))\cup\{\infty\}$.
Let $da$ be the one-form dual to the fundamental vector field $H^*$ on $A$
corresponding to
$H\in \aaaa^+$ with $|H|=1$.
We decompose forms $\omega$ in $\cL^p$ as $\omega=\omega_1\oplus da\wedge
\omega_2$,
where $i_{H^*}\omega_j=0$, $j=1,2$, and $i_{H^*}$ is the insertion of $H^*$.
One can check that
$$\underline{S_\infty\cL^p[G]}(W_i)=V^p(W_i)\oplus da\wedge V^{p-1}(W_i)\ ,$$
where
\begin{eqnarray*}
V^p(W_i)&:=&\{\omega\in C^\infty(U_i\times (i,\infty)\times
G,\pi_1^*\Lambda^pT^*U_i)|\\
&&\forall L_1\in\cU(\naaa),L_2\in\cU(\aaaa),L_3\in\cU(\gaaa)\:\:\:\exists
R,Q\in\R\\
&&s.t.\:\:\:\sup_{a\in(i,\infty),h\in G}  a^{-R} \ee^{-Q\dist(\orig,hK)}
|\omega(L_1n,L_2a,L_3h)|<\infty\}\ ,
\end{eqnarray*}
and $\pi_1:U_i\times (i,\infty)\times G\rightarrow U_i$ is the projection onto
the first factor.
The complex $(\underline{S_\infty\cL^*[G]}(W_i),d)$
is the total complex
of a double complex. The latter consists of  two rows equal
$(V^*(W_i),d)$, where $d$ is the differential of the de Rham complex
of $U_i$. The vertical differential of the double complex
is the differentiation along the $A$-direction given by $\pm H$.
The balls $U_i$ can be contracted radially.
The associated zero homotopy of the de Rham complex of $U_i$
extends to $V^*(W_i)$.
Again we leave the verification to the reader (see also the proof of Lemma
\ref{amdt}).
Thus the rows $(V^*(W_i),d)$ of the double complex are acyclic.
The zeroth horizontal cohomology is equal to ${}^\naaa V^0(W_i)$,
the functions in $V^0(W_i)$ that do not depend on $n\in U_i$.
To finish the proof of the Lemma we must show that
$H:{}^\naaa V^0(W_i)\rightarrow {}^\naaa V^0(W_i)$
is surjective.
In fact the equation $Hf=g$, $g\in {}^\naaa V^0(W_i)$, can explicitly be
solved by integration
such that $f\in {}^\naaa V^0(W_i)$. Set $$f(n,a,h)=\int_i^a g(n,b,h) db\ .$$
\hB
 Let $\pi:X\times \Gamma\backslash G$ be the projection onto the first factor.
We define
\begin{eqnarray*}
S_\infty\cL^p[\Gamma\backslash G]&:=&\{\omega\in C^\infty(X\times
\Gamma\backslash  G ,\pi^*L^*)\\
&&|\:\forall L_1,L_2\in \cU(\gaaa)\:\:\: \exists R,Q\in\R\\
&&s.t.\:\: \sup_{g,h \in G}\ee^{-R\dist(\orig,gK)}
\ee^{-Q\dist_Y(\Gamma\orig,\Gamma hK)} |\omega(L_1g,L_2h)|<\infty
\} \ ,
\end{eqnarray*}
where $\dist_Y(\Gamma\orig,\Gamma hK)$ is the distance of $\Gamma hK$ from
$\Gamma \orig$ in $Y:=\Gamma\backslash X$.
\begin{lem}\label{mglobal}
The de Rham complex
$(S_\infty\cL^*[\Gamma\backslash  G ],d)$ is acyclic.
\end{lem}
\proof
As in the proof of Lemma \ref{mmarrt} the de Rham complex
$(S_\infty\cL^*[\Gamma\backslash  G ],d)$ is the total complex of a double
complex
consisting of two rows each equal to the $\naaa$-cohomology complex
$(\Lambda^*\naaa\otimes S_\infty\cL^0[\Gamma\backslash  G ],d)$.
Let $\pi_1:N\times A\times \Gamma\backslash  G \rightarrow N$ be the projection
onto the first
factor.
We identify
$\Lambda^p\naaa\otimes S_\infty\cL^0[\Gamma\backslash  G ]$ with
\begin{eqnarray*}
 &&\{\omega\in C^\infty(N\times A\times \Gamma\backslash  G
,\pi_1^*\Lambda^pT^\ast N)
\\&&|\:\: \forall  L_1\in\cU(\naaa), L_2\in \cU(\aaaa),L_3\in\cU(\gaaa)\:\:\:
\exists R,Q\in \R\\
&& s.t. \sup_{na\in NA,g\in    G , I\in I_p}
\ee^{-R\dist(\orig,naK)}\ee^{-Q\dist_Y(\Gamma \orig,\Gamma gK)}
|\omega_I(L_1n,L_2a,L_3g)| < \infty\}\ .
\end{eqnarray*}
We show that the complex $(\Lambda^*\naaa\otimes S_\infty\cL^0[\Gamma\backslash
 G ],d)$ is acyclic.
In the proof of Lemma \ref{amdt} we constructed a zero homotopy $H=\int_0^1
H_t$, where
$(H_t\omega)(n,a,g)=\Psi_t^\ast(i_{T_t}\omega(\Psi_t(n),a,g))$.

We claim that $H:\Lambda^p\naaa\otimes S_\infty\cL^0[\Gamma\backslash  G
]\rightarrow \Lambda^{p-1}\naaa\otimes S_\infty\cL^0[\Gamma\backslash  G ]$.
The same discussion as in the proof of Lemma \ref{amdt} leads to the following
estimate.
Fix $L_1\in\cU(\naaa)$, $L_2\in\cU(\aaaa)$, $L_3\in\cU(\gaaa)$.
Let $r:=\deg(L_1)$.
For all
$\omega$ with
$$\sup_{l\in L_r,na\in NA,g\in G , I\in I_p} \ee^{-P\dist_Y(\Gamma \orig,\Gamma
gK)}\ee^{-R\dist(\orig,naK)} |\omega_I(X(l)n,L_2a,L_3g)|=:S(\omega) < \infty$$
we have
$$ \sup_{na\in NA,g\in  G ,J\in I_{p-1}}
\ee^{-P\dist_Y(\Gamma \orig,\Gamma gK)}\ee^{-(R+Q)\dist(\orig,naK)}
|(H_t\omega)_J(L_1n,L_2a,L_3g)| < CS(\omega)$$
with $C<\infty$ independent of $\omega$.
This easily implies the claim.
Hence the rows of the double complex are acyclic
and their zeroth cohomology is  equal to
\begin{eqnarray*}
V&:=&\{f\in C^\infty(A\times \Gamma\backslash  G )\:\: |\:\: \forall  L_2\in
\cU(\aaaa), L_3\in\cU(\gaaa)\:\:\: \exists R,Q\in \R\\
&&\hspace{3cm} s.t. \sup_{a\in A,g\in G } \ee^{-R|\log(a)|}
\ee^{-Q\dist_Y(\Gamma \orig,\Gamma gK)} |f(L_2a,L_3g)| < \infty\}\ .
\end{eqnarray*}
The vertical differential given by $H:V\rightarrow V$ is  surjective.
In fact, let $f\in V$ and set
$$F(a,g):=\int_1^a f(b,g)db\ .$$
Then $HF=f$ and $F\in V$.
This proves the lemma. \hB

 \newcommand{\cHC}{{{\cal HC}(\gaaa,K)}}

\section{The standard resolution}
For the convenience of the reader we repeat here the construction
of the standard resolution given in \cite{bunkeolbrich947}.

Let $(\pi,V_{\pi,K})\in\cHC$ be a Harish-Chandra module.
Then $V_{\pi,K}$ decomposes into a direct sum of joint
generalized eigenspaces of $\cZ(\gaaa)$.
Since the summands can be treated separately, without loss of
 generality we can assume
that there exist $\mu\in \C$ and $k\in\Naaa$ such
that $B:=(\Omega-\mu)^k\in \Ann(V_{\pi,K})$, i.e., $BV_{\pi,K}=0$.

Let $W$ be a finite-dimensional $K$-stable subspace of
 the dual $V_{\tilde\pi,K}$  of $V_{\pi,K}$ in the category
$\cHC$, which generates $V_{\tilde\pi,K}$ as a
$\cU(\gaaa)$-module. Let $E_0\rightarrow X$ be
the homogeneous vector bundle $G\times_K \tilde W$
and $\cE_0$ be the space of its smooth sections.
Using any globalization $ V_\pi $ of $V_{\pi,K}$ (i.e. a
representation of $G$ such that $ V_\pi=V_{\pi,K}$)
we can define an embedding
$$
i: V_{\pi,K}\hookrightarrow \cE_0\cong [C^\infty(G)\otimes \tilde W]^K
$$
that is characterized by
$$\langle i(v)(g),w\rangle:= \langle w,\pi(g^{-1})v\rangle,
\qquad v\in V_{\pi,K},w\in W,g\in G.$$
In fact $i$ maps into $S_\infty\cE_0$
and the closure of $i(V_{\pi,K})$ in $S_\infty\cE_0$ is contained in
$S_\infty\cE_0(B)$ and constitutes the distribution vector
globalization $V_{\pi,-\infty}$
of $V_{\pi,K}$ (Wallach \cite{wallach92}, Ch. 11, Casselman
\cite{casselman89}).

We will also consider the space $V_{\pi,for}:=V_{\tilde\pi, K}^*$
of formal power series vectors
of $V_{\pi,K}$. There is an exact functor from $\cHC$ to the
category of (not necessarily $K$-finite) $(\gaaa,K)$-modules
which sends $V_{\pi,K}$ to $V_{\pi,for}$. Note that
$V_{\pi,for}=\prod_{\gamma\in\hat K} V_{\pi,K}(\gamma)$.

For homogeneous vector bundles $E$ and $F$ on $X$ we denote by $D(E,F)$ the set
of $G$-invariant differential operators $D:\cE\rightarrow \cF$.
\begin{prop}\label{1.zeile}
There exist homogeneous vector bundles $E_1,E_2,\dots$ on $X$ and
$G$-\linebreak[4] invariant differential operators $D_i\in D(E_i,E_{i+1})$,
$i=0,1,\dots$, such that the embedding $i:V_{\pi,-\omega}\hookrightarrow
\cE_0(B)$ can be
extended to a (possibly infinite) exact sequence
\begin{equation}\label{*}
 0\rightarrow V_{\pi,-\infty}\stackrel{i}{\rightarrow}S_\infty\cE_0(B)
\stackrel{D_0}{\rightarrow}S_\infty\cE_1(B)
\stackrel{D_1}{\rightarrow}S_\infty\cE_2(B)
\stackrel{D_2}
{\rightarrow}\dots\ .
\end{equation}
This sequence remains to be exact on the level of formal power
series :
\begin{equation}\label{*+}
 0\rightarrow V_{\pi,for}\stackrel{i}{\rightarrow}\cE_0^{for}(B)
\stackrel{D_0}{\rightarrow}\cE_1^{for}(B)
\stackrel{D_1}{\rightarrow}\cE_2^{for}(B)
\stackrel{D_2}{\rightarrow}\dots\ .
\end{equation}
\end{prop}
\proof
Let $\cZ(E)$ be the image of $\cZ(\gaaa)$ in $D(E,E)$.
For any vector bundle $E\rightarrow X$ the $\C[B]$-module $\cZ(E)$
is finitely generated (\cite{bunkeolbrich947}, Lemma 2.3).
\begin{lem}\label{isntit}
For any  vector bundle $E\rightarrow X$ we have $\cE(B)_K\in\cHC$.
\end{lem}
\proof
Let $(\gamma,V_\gamma)$ be the finite dimensional representation
of $K$ corresponding to $E$ and $(\tilde{\gamma},V_{\tilde{\gamma}})$
be  its dual.
We consider the $K$-equivariant embedding
$$i:V_{\tilde{\gamma}}\hookrightarrow \widetilde{\cE(B)_K}$$
defined by
$$i(\tilde{v})(f):=\langle\tilde{v},f(e)\rangle\ ,$$
where we identify the fibre of $E$ at $e=[K]$
with $V_\gamma$. Let $T:=\cU(\gaaa)(i(V_{\tilde{\gamma}}))$.
For any $t\in T$  the dimension of $\cZ(\gaaa)t$ can be estimated
by the dimension of a generating subspace of the
 $\C[B]$-module $\cZ(E)$. Thus $T$ is a locally
$Z(\gaaa)$-finite and finitely generated $\cU(\gaaa)$-module. By a theorem of
Harish-Chandra (\cite{wallach88}, 3.4.7), $T\in\cHC$. The canonical map
$\cE(B)_K\rightarrow \tilde{T}$
is injective by the analyticity of solutions of the equation $Bf=0$. In fact,
an element in the kernel of this map
would have a vanishing Taylor series at $e$.
We obtain that $T\hookrightarrow \widetilde{\cE(B)_K}$ is surjective.
Thus $T=\widetilde{\cE(B)_K}$ and $\cE(B)_K\in\cHC$ since
the dual of a Harish-Chandra module is a Harish-Chandra module, too
(\cite{wallach88}, 4.3.2).
\hB
\begin{lem}\label{hno}
Let $V_{\pi,K}$ be a Harish-Chandra submodule of $\cE(B)_K$.
Then there exist a homogeneous vector bundle $F$ and an
 operator $D\in D(E,F)$ such that $\ker D \cap \cE(B)_K=
V_{\pi,K}$. We also have $\ker D \cap S_\infty\cE(B)=V_{\pi,-\infty}$.
\end{lem}
\proof
According to the proof of Lemma \ref{isntit}
there is a surjection
$$ \cU(\gaaa)\otimes_{\cU(\kaaa)} V_{\tilde\gamma}
\rightarrow \widetilde{\cE(B)_K}\ .$$ Let $W$ be a finite-dimensional
$K$-stable generating subspace of the Harish-Chandra module
$V_{\pi,K}^\perp\subset\widetilde{\cE(B)_K}$. Then we choose a
$K$-equivariant map $\alpha$ such that the following diagram
$$\begin{array}{ccc}
W&\stackrel{\alpha}{\longrightarrow}&\cU(\gaaa) \otimes_{\cU(\kaaa)}
V_{\tilde\gamma}\\
\downarrow&&\downarrow\\
V_{\pi,K}^\perp&\longrightarrow &\widetilde{\cE(B)_K}
\end{array}$$
commutes. This is possible since
$\cU(\gaaa)\otimes_{\cU(\kaaa)}V_{\tilde\gamma}$ is
$K$-semisimple.

We set $F:=G\times_K \tilde W$. The map $\alpha$ can be considered
as an element of
$$[\cU(\gaaa)\otimes_{\cU(\kaaa)} V_{\tilde\gamma}\otimes
\tilde W]^K\cong [\cU(\gaaa)\otimes_{\cU(\kaaa)}
\Hom (V_{\gamma},\tilde W)]^K\ .$$
The latter space is canonically isomorphic to $D(E,F)$ via the
right regular  representation $R$ of $\cU(\gaaa)$ on $C^\infty(G)\otimes
V_\gamma$. Thus $\alpha$ defines an element
$D\in D(E,F)$. If $\alpha(w)=\sum X_i\otimes v_i$, then
$$ \langle w,Df\rangle_F = \sum \langle v_i, R_{X_i}f\rangle_E \in
C^\infty(G),\ w\in W,v_i\in V_{\tilde\gamma}, X_i\in \cU(\gaaa)\ .$$
Let $f\in \cE(B)$, $X\in\cU(\gaaa)$ and $w\in W$.
Then we have
\begin{eqnarray}\label{mufti}
\langle w, L_XDf(1)\rangle_F&=&\langle w, DL_Xf\rangle_F\nonumber\\
&=&\sum\langle v_i,R_{X_i}L_Xf(1)\rangle_E\\
&=&\langle w ,L_Xf \rangle_{\cE(B)}\nonumber\\
&=& \langle L_{X^{op}}w ,f \rangle_{\cE(B)}\ ,\nonumber
\end{eqnarray}
where $X\rightarrow X^{op}$ is the anti-automorphism of
$\cU(\gaaa)$ induced by the multiplication with $-1$ on
$\gaaa$. By construction  $Df=0$ iff the left hand side of (\ref{mufti})
vanishes for all $X\in\cU(\gaaa)$ and $w\in W$,
while $f\in V_{\pi,-\infty}$ iff the right hand side
does. The lemma follows.\hB

In order to construct the bundles
$E_i$ and operators $D_i$ of Proposition \ref{1.zeile} we
iterate Lemma \ref{hno}. $D_i(\cE_i(B)_K)$ is a Harish-Chandra submodule of
$\cE_{i+1}(B)_K$.
Therefore we find a bundle $E_{i+2}$ and an operator $D_{i+1}\in
D(E_{i+1},E_{i+2})$ such that $\ker D_{i+1}\cap\cE_{i+1}(B)_K =
D_i(\cE_i(B)_K)$. We obtain an exact sequence of Harish-Chandra modules
$$ 0\rightarrow V_{\pi,K}\stackrel{i}{\rightarrow}\cE_0(B)_K
\stackrel{D_0}{\rightarrow}\cE_1(B)_K  \stackrel{D_1}{\rightarrow}\cE_2(B)_K
\stackrel{D_2}{\rightarrow}\dots\ .$$
Applying the distribution vector globalization functor (which
is exact) we end up with (\ref{*}). Analogously, we want to obtain (\ref{*+})
by taking
formal power series vectors. This is possible since for any homogeneous vector
bundle $E$ we have
$$ (\cE(B)_K)_{for}=\cE^{for}(B)$$
(\cite{bunkeolbrich947}, Lemma 3.5).

The Proposition \ref{1.zeile} provides a resolution
of $V_{\pi,-\infty}$ by spaces $S_\infty\cE_i(B)$.
We now employ a Koszul complex construction in order to get rid of
the eigenspaces. We recall the following fact from \cite{bunkeolbrich947}.
\begin{lem}\label{mlem}
Let $E$, $F$ be homogeneous vector bundles on $X$ and
$A\in D(E,F)$
such that $A  \cE(B) =0$.
Then $A = H B$ for some $H\in D(E,F)$.
\end{lem}
Let $V_{\pi,K},E_i,D_i$ be constructed as in Proposition \ref{1.zeile}.
\begin{prop}\label{stare}
There exist $H_i\in D(E_i,E_{i+2})$, $i\ge 0$, making the
following into an exact complex:
\begin{equation}\label{**}
0\rightarrow V_{\pi,-\infty}\rightarrow S_\infty\cE_0
\stackrel{{\scriptsize \left(\begin{array}{c}D_0\\B\end{array}\right)}}
{\longrightarrow}
\begin{array}{c}S_\infty\cE_1\\ \oplus\\ S_\infty\cE_0\end{array}
\stackrel{ \left({\scriptsize\begin{array}{cc}D_1&H_0\\-B&D_0\end{array}}
\right)}{\longrightarrow}
\begin{array}{c}S_\infty\cE_2\\ \oplus\\S_\infty\cE_1\end{array}
\stackrel{ \left({\scriptsize\begin{array}{cc}
D_2&H_1\\B&D_1\end{array}}\right)}{\longrightarrow}\dots\ .
\end{equation}
\end{prop}
We shall call (\ref{**}) a standard resolution of $V_{\pi,-\infty}$.\\
\proof
In order to construct the operators $H_i$ we apply
Lemma \ref{mlem} for $A=D_{i+1}D_i$.
The exactness of (\ref{**}) is easily reduced to the exactness of
 (\ref{*}) and the surjectivity
of $B:S_\infty\cE_i\rightarrow S_\infty\cE_i$ proved in Theorem \ref{sur}.
\hB

    \section{$\naaa$-cohomology
% of globalizations of Harish-Chandra modules
}
 Let $B=(\Omega_G-\lambda)^l$ for some $\lambda\in\C$, $l\in\nat$, where
$\Omega_G$
is the Casimir operator of $G$,
and $S_\infty\cE(B)=\{f\in S_\infty\cE\:|\: Bf=0\}$.
\begin{lem}\label{lez}
We have
$$H^p(\naaa,S_\infty\cE(B))=0,\quad \forall p\ge 1.$$
\end{lem}
\proof
By    Theorem \ref{sur} and Lemma \ref{amdt}
$$0 \rightarrow S_\infty\cE(B)\rightarrow S_\infty
\cE\stackrel{B}{\rightarrow}S_\infty \cE\rightarrow 0$$
is an $\naaa$-acyclic resolution of $S_\infty\cE(B)$.
Taking $\naaa$-invariants we obtain the complex
\begin{equation}\label{ueeqq1}0\rightarrow {}^\naaa S_\infty \cE(B)\rightarrow
S_\infty  C^\infty(A)\otimes V_\gamma\stackrel{{}^\naaa B}{\rightarrow}
S_\infty C^\infty(A)\otimes V_\gamma\rightarrow 0\ .\end{equation}
Here ${}^\naaa B$ is the restriction of $B$
to the subspace of $\naaa$-invariant vectors.
It is a second order translation invariant differential operator on $A$.
The complex (\ref{ueeqq1}) is again exact and the Lemma follows.
\hB

Let $(\pi,V_{\pi,K}) \in\hc$.
Recall that $H^*(\naaa,V_{\pi,-\infty})$ carries a natural $MA$-module
structure.
 \begin{theorem}\label{haha}
The inclusion $V_{\pi,-\infty}\hookrightarrow V_{\pi,for}$
induces an isomorphism
$$H^p(\naaa,V_{\pi,-\infty})\stackrel{\sim}{\longrightarrow}
H^p(\naaa,V_{\pi,for})\ .$$
Moreover
$H^p(\naaa,V_{\pi,-\infty})=H^p(\naaa,V_{\pi,-\omega})= H^p(\naaa,V_{\pi,for})$
and all spaces are finite dimensional.
The $\naaa$-cohomology of $V_{\pi,-\infty}$ satisfies Poincar\'e duality
\begin{equation}\label{ukk1}
H^p(\naaa,V_{\pi,-\infty})^*\cong H^{\dim(\naaa)-p}(\naaa,
V_{\tilde{\pi},\infty})\otimes\Lambda^{\dim{\naaa}}\naaa .
\end{equation}
We also have $$ H^p(\naaa,V_{\pi,\infty})\cong H^p(\naaa,V_{\pi,K}).
$$
 \end{theorem}
\proof
By (\cite{bunkeolbrich947}, Lemma 2.3 and Proposition 4.1), Lemma \ref{amdt}
and Proposition \ref{1.zeile} the cohomology  $H^p(\naaa,V_{\pi,*})$ for
$*=-\infty$, $for$ is isomorphic to the cohomology of the subcomplex of
$\naaa$-invariants of (\ref{*}), (\ref{*+}), respectively.
Hence the following
lemma implies the first assertion of the theorem.
\begin{lem}\label{hihi}
For any homogeneous vector bundle $E\rightarrow X$ associated to $V_\gamma$ we
have
$$ {}^\naaa S_\infty \cE(B)={}^\naaa\cE^{for}(B) \ .$$
\end{lem}
\proof
The $\cU(\aaaa)$-module
$${}^\naaa \cE^{for}(B)\cong (\widetilde {\cE(B)_K}/\naaa(\widetilde
{\cE(B)_K}))^*$$
is finite dimensional (see \cite{wallach88}, Ch.4). Therefore it splits into
generalized weight spaces ${}^\naaa \cE^{for}(B)_\mu, \mu\in\aaaa^*_\C$. $f\in
{}^\naaa \cE^{for}(B)_\mu$,
considered as a formal power series on $\aaaa$,
satisfies the differential equation
\begin{equation}\label{muuh}(H+\mu(H))^kf=0\quad  \forall
H\in\aaaa\end{equation}
for a certain $k\in \nat$.
The solutions of (\ref{muuh}) have the form
$$P(H)\ee^{-\mu(H)},\qquad P\in S(\aaaa^*)\otimes V_\gamma\ .$$
They extend to smooth $\naaa$-invariant sections in ${}^\naaa S_\infty \cE(B)$.
\hB
The proof of the remaining assertions of the theorem is parallel
to the proofs of the corresponding facts in \cite{bunkeolbrich947}.
\hB

  \section{$S_\infty\cE$ is $\Gamma$-acyclic}
Let $\Gamma\subset G$ be a discrete, torsion-free subgroup of finite covolume.
Let $E\rightarrow X$ be a $G$-homogeneous vector bundle
and $S_\infty\cE$ the space of its sections of moderate growth.
\begin{theorem}\label{gaz}
$S_\infty\cE$ is  $\Gamma$-acyclic, i.e.,
$$H^p(\Gamma,S_\infty\cE)=0\quad\forall p\ge 1\ .$$
\end{theorem}
\proof
We first consider the space of functions of moderate growth
$S_\infty C^\infty(G)$ on $G$ defined by
$$S_\infty C^\infty(G):=\{f\in\C^\infty(G)\:|\: \forall L\in\cU(\gaaa)\:\exists
R\in\R\:s.t. \:
\sup_{g\in G} \ee^{-R\dist(gK,\orig)} |f(Lg)| <\infty\}\ .$$
As a topological vector space $S_\infty C^\infty(G)$ is a limit of Fr\'echet
spaces.
\begin{prop}\label{acc}
$S_\infty C^\infty(G)$ is $\Gamma$-acyclic.
\end{prop}
\proof
$H^\ast(\Gamma,S_\infty C^\infty(G))$ is the cohomology of the
de Rham complex
of $Y=\Gamma\backslash X$ twisted with
the flat bundle associated to the $\Gamma$-module
$S_\infty C^\infty(G)$.
In greater detail let $L^\ast:=\Lambda^\ast T^\ast X$ and
 $\cC^\ast:=C^\infty(X,L^\ast\otimes S_\infty C^\infty(G))$.
Moreover, let $d:\cC^\ast\rightarrow \cC^{\ast+1}$ denote the  differential of
the de Rham complex.
The complex $(\cC^\ast,d)$ is a complex of $\Gamma$-modules.
If we view $\omega\in\cC^p$ as a function on $G$ with values in $\cL^p$,
then the action of $\gamma\in\Gamma$ on $\omega$ is given by
$(\gamma\omega)(g)=(L_\gamma^*\omega)(\gamma^{-1}g)$,
where $L_\gamma^\ast$ is the pull back of forms associated to the
diffeomorphism $L_\gamma:X\rightarrow X$
given by $L_\gamma(x):=\gamma^{-1}x$.
The complex $(\cC^\ast,d)$ is exact
(the contraction of $X$ along radial rays induces a zero-homotopy
of the de Rham complex)
and it consists of $\Gamma$-acyclic modules (\cite{bunkeolbrich947}, Lemma
2.4).
Hence $H^\ast(\Gamma,S_\infty C^\infty(G))$ is the cohomology of the complex
$({}^\Gamma\cC^\ast,d)$.

Let $\cS^\ast:=S_\infty \cL^\ast[S_\infty(G)]$
(see Section \ref{sec15} for notation).
Then $(\cS^\ast,d)\hookrightarrow (\cC^\ast,d)$
is a subcomplex.
\begin{lem}\label{incl}
The inclusion $({}^\Gamma\cS^\ast,d)\hookrightarrow ({}^\Gamma\cC^\ast,d)$
induces an isomorphism in cohomology.
\end{lem}
\proof
The manifold $Y$ has finitely many cusps each diffeomorphic to
$B_i\times [0,\infty)$, where $B_i$ is some compact nil-manifold.
The Borel-Serre compactification of $Y$  is obtained
by attaching copies of the cusp bases $B_i$
as a boundary at infinity obtaining a manifold with boundary $\bar{Y}$.
There is a natural sheafification
$(\underline{{}^\Gamma\cS^\ast},d)\hookrightarrow
(\underline{{}^\Gamma\cC^\ast},d)$
of the inclusion $({}^\Gamma\cS^\ast,d)\hookrightarrow ({}^\Gamma\cC^\ast,d)$
on $\bar{Y}$.
To any finite open covering of $\bar{Y}$ there is an associated
partition of unity
which is compatible with the sheafs $\underline{{}^\Gamma\cS^\ast}$.
Thus $\underline{{}^\Gamma\cS^\ast}$ and $\underline{{}^\Gamma\cC^\ast}$
are acyclic with respect to the global section functor.

The complex of sheaves corresponding to $({}^\Gamma\cC^\ast,d)$
is locally acyclic by the standard Poincar\'e Lemma.
Moreover the inclusion of sheaves
$(\underline{{}^\Gamma\cS^\ast},d)\hookrightarrow
(\underline{{}^\Gamma\cC^\ast},d)$
induces an isomorphism of the zeroth
cohomology sheaves.
By Lemma  \ref{mmarrt} $(\underline{{}^\Gamma\cS^\ast},d)$
is locally acyclic.
Thus the inclusion
$(\underline{{}^\Gamma\cS^\ast},d)\hookrightarrow
(\underline{{}^\Gamma\cC^\ast},d)$
is a quasi-isomorphism.
Since the sheaves $\underline{{}^\Gamma\cS^\ast},\underline{{}^\Gamma\cC^\ast}$
are acyclic with respect
to the global section functor, the induced map of the complexes
of  global sections $({}^\Gamma \cS^*,d)\hookrightarrow ({}^\Gamma \cC^*,d)$ is
a quasi-isomorphism, too.
\hB
\begin{lem} \label{moac}
$H^p({}^\Gamma\cS^.,d)=0$, $\forall p\ge 1$.
\end{lem}
\proof
The map $T:X\times G\rightarrow X\times G$ given by
$(x,g)\rightarrow (gx,g)$ intertwines the $\Gamma$-action on the second factor
$G$ with the diagonal $\Gamma$-action
on the product $X\times G$.
 We claim that $T$ induces an isomorphism
$T^\ast:S_\infty\cL^\ast[S_\infty(G)]\rightarrow S_\infty\cL^\ast[S_\infty(G)]$
intertwining the  $\Gamma$-module structure given above with the
$\Gamma$-module structure
given by $(\gamma\omega)(g)=\omega(\gamma^{-1}g)$
(again viewing $\omega$   as a function from $G$ to $\cL^*$).
 The inverse of $T^{\ast}$ is induced by $T^{-1}:X\times G\rightarrow X\times
G$,
  $T^{-1}(x,g)=(g^{-1}x,g)$.
In fact, $(T^\ast\omega)(g)=L_{g^{-1}}^\ast\omega(g)$ and
$$(T^\ast\gamma\omega)(g)=L_{g^{-1}}^\ast L_\gamma^\ast
\omega(\gamma^{-1}g))=(T^\ast\omega)(\gamma^{-1}g)\ . $$
In order to prove the claim we must show that $T^\ast,(T^*)^{-1}$
are compatible with the weighted spaces.
 It is at this point that we have  to consider the weighted
de Rham complex $\cS^*$.
Since $T$ mixes the $G$- and the $X$ directions
it does not act on $\cC^\ast$.

Let $\Delta:\cU(\gaaa)\rightarrow \cU(\gaaa)\otimes \cU(\gaaa)$
be the co-product induced by $X\mapsto X\otimes 1 + 1 \otimes X$, $X\in \gaaa$.
Fix  $L_1,L\in \cU(\gaaa)$ and let $(\Delta\otimes \id )\Delta(L) =\sum_\alpha
A_\alpha\otimes B_\alpha\otimes C_\alpha$.
Let $\omega\in \cS^p$.
Then $$(T^\ast\omega)(L_1h,Lg)=\sum_{\alpha}  \omega(B_\alpha gL_1h,C_\alpha
g)\circ DL_{g^{-1}A^{op}_\alpha}\ ,$$
where $DL_{g}$ is the differential of $L_g$ acting on $TX$.

 Let $r=\deg(L_1)+\deg(L)$.  Assume that
$$\sup_{g,h\in G, l,l_1\in L_r,}
\ee^{-R\dist(\orig,gK)}\ee^{-Q\dist(\orig,hK)}  |\omega(X(l)h,X(l_1)g)|<\infty
$$
for appropriate $Q,R\in \R$.
Here $\{X(l)\}_{l\in L_r}$ is a basis of the differential operators on $G$ of
order $\le r$
 (similarly to the notation in the proof of Lemma \ref{amdt}).
Note that $\dist(\orig,ghK)\le \dist(\orig,hK)+ \dist(\orig,gK)$.
Moreover we have
$$|DL_{g^{-1}A^{op}_\alpha}|\le C \ee^{-P\dist(\orig,gK)} \quad \forall \alpha,
\:\:\forall g\in G $$
and
$B_\alpha gL_1h=\sum_{l\in L_r} e_{\alpha,l}(g) X(l)gh$ with
$$|e_{\alpha,l}(g)|\le C \ee^{-P\dist(\orig,gK)}
 \quad \forall l\in L(r),\:\:\:\forall \alpha, \:\:\forall g\in G$$
for $P\in \R$ large enough.
We conclude that
$$\sup_{g,h\in G} |(T^\ast\omega)(L_1h,Lg)|
\ee^{-(R+Q+2P)\dist(\orig,gK)}\ee^{- Q \dist(\orig,hK)}
<\infty\ .$$
Hence $T^*\omega\in \cS^p$.
In a similar manner one can handle $(T^*)^{-1}$ thus proving the claim.

We see that  $H^\ast({}^\Gamma \cS^.,d)$ is isomorphic to the cohomology of
$(S_\infty\cL^\ast[  \Gamma \backslash G],d)$.
Lemma \ref{moac} now follows from Lemma \ref{mglobal}. \hB

We now finish the proof of Theorem \ref{gaz}.
Note that $S_\infty C^\infty(G)$ carries a right $K$-module structure which
commutes with the left $\Gamma$-module structure.
This induces a $K$-action on $({}^\Gamma \cC^*,d)$.
$H^*(\Gamma,S_\infty\cE)$ is the cohomology of the complex
$([{}^\Gamma \cC^\ast\otimes V_\gamma]^K,d)$.

Let $[z]\in H^p(\Gamma, S_\infty \cE)$, $p\ge 1$,
be represented by the $K$-invariant cycle $z\in [{}^\Gamma \cC^p\otimes
V_\gamma]^K$.
Then by Lemma  \ref{acc} we have
$z=d b$ for some possibly non-invariant $p-1$ cochain
$b\in  {}^\Gamma \cC^{p-1}\otimes V_\gamma$. Let $\bar{b}$ the
average of $b$ with respect to $K$. Then also $d\bar{b}=z$ and hence
$0=[z]\in H^p(\Gamma, S_\infty\cE)$.
This proves Theorem  \ref{gaz} .\hB
  \section{The cokernel of $B:S_\infty\cE_Y\rightarrow
S_\infty\cE_Y$}\label{klop}
In this section we relate the dimension of the cokernel of $B:=\Omega-\mu$ on
$S_\infty\cE_Y:={}^\Gamma S_\infty\cE$ with the
dimension of corresponding spaces of cusp forms.
Since it is not a priori clear that the range of $B$ is closed note that
we consider the algebraic cokernel of $B$.

Let $E_Y:=\Gamma\backslash E$ be the locally homogeneous
bundle over $Y$. A section $f\in C^\infty(Y,E_Y)=:\cE_Y$
can be viewed as a function on $f:G\rightarrow V_\gamma$
satisfying $f(h g  k)=\gamma^{-1}(k)f(g)$ for all $h\in \Gamma$, $k\in K$.
The space $S_\infty\cE_Y$ has the alternative description
$$
S_\infty\cE_Y=\{f\in \cE_Y\:|\: \forall L\in\cU(\gaaa)\: \exists R\in \R \:\:
s.t. \:\:\sup_{g\in G}\ee^{R\dist_Y(\Gamma\orig,\Gamma gK)}|f(Lg)| <\infty\}\ .
$$
A section $f\in \cE_Y$ which is an  eigenfunction of
$\Omega$ and satisfies
$$\sup_{g\in G}   \ee^{R\dist_Y(\Gamma\orig,\Gamma gK)}
|f(g)|<\infty\:\:\:\forall R\in\R$$
is called a cusp form.
Let $V_\mu$ be the space of cusp forms in $\ker(B)$.
If $V_\mu\not=\{0\}$, then  necessarily $\mu\in\ R$, since $\Omega$
is symmetric.

The main result of the present section is the following theorem.
\begin{theorem}\label{koko}
$$\dim \coker(B:S_\infty\cE_Y\rightarrow S_\infty\cE_Y)=\dim V_\mu\ .$$
 \end{theorem}
The proof of the theorem occupies the remainder of the section.
\begin{lem}\label{larger}
 $$\dim \coker(B:S_\infty\cE_Y\rightarrow S_\infty\cE_Y)\ge \dim V_\mu\ .$$
 \end{lem}
\proof
We show that the projection of $V_\mu$ to $\coker(B:S_\infty\cE_Y\rightarrow
S_\infty\cE_Y)$
is an inclusion.
We can assume that $V_\mu\not=0$ and hence $\mu\in\R$.
Thus $B$ is symmetric.

If $f\in V_\mu$, then it vanishes rapidly and can be integrated against
elements
of $S_\infty\cE_Y$.
Thus let $f\in V_\mu$ and assume that $[f]=0$ in
$\coker(B:S_\infty\cE_Y\rightarrow S_\infty\cE_Y)$.
Then $f=Bh$ for some $h\in S_\infty\cE_Y$.
We have
\begin{eqnarray*}
0&=&\int_Y  \langle (Bf)(y)\:,\:h(y)\rangle_{E_{Y,y}} dy \\
&=&\int_Y \langle f(y)\:,\:(Bh)(y)\rangle_{E_{Y,y}} dy\\
&=&  \|f\|^2_{L^2(Y,E_Y)}\ .
\end{eqnarray*}
It follows $f=0$.
This proves Lemma \ref{larger}. \hB
We define an increasing sequence of Fr\'echet spaces $S_R\cE_Y$, $R\in\R$, by
$$S_R\cE_Y:=\{f\in\cE_Y\:|\: p_{Y;-R,L}(f)<\infty\: \:\forall L\in\cU(\gaaa)\}\
,$$
where the seminorms are defined by
$$p_{Y;R,L}(f):=\sup_{g\in G} \ee^{R\dist(\Gamma\orig,\Gamma gK)} |f(gL)|\ .$$
Note that (in contrast to the definition of $S_\infty\cE$)
we employ  the right action of $\cU(\gaaa)$ to define $p_{Y;R,L}(f)$.
In fact, to define $S_R\cE_Y$ it is sufficient to consider the seminorms
$p_{Y;R,\Omega^k}(f)$,
$k\in\nat_0$.
Let
$$
S_{-\infty}\cE_Y = \bigcap_{R\in\R}S_R\cE_Y$$
  with the natural topology
of the intersection and $S_\infty\cEp_Y$  be the topological conjugate  dual of
$S_{-\infty}\cE_Y$.
The Hermitian scalar product of $E_Y$ induces an inclusion
$S_\infty\cE_Y\hookrightarrow S_\infty\cEp_Y$.
 \begin{lem}\label{smmo}
The inclusion $S_\infty\cE_Y\hookrightarrow S_\infty\cEp_Y$
induces an injection
$$\coker(B:S_\infty\cE_Y\rightarrow S_\infty\cE_Y)\hookrightarrow
\coker(B:S_\infty\cEp_Y\rightarrow S_\infty\cEp_Y)\ .$$
\end{lem}
\proof
We reduce the proof to Lemma \ref{reg}.
Let $\pi:X\rightarrow Y$ be the projection.
We define a continuous map $\pi_*:S_{-\infty}\cE\rightarrow S_{-\infty}\cE_Y$.
For $f\in S_{-\infty}\cE$ let $(\pi_\ast f)(\Gamma g):=\sum_{\gamma\in \Gamma}
f(\gamma g)$.
We show that the sum converges and $\pi_\ast$ is continuous.

There is a $Q\in R$ such that
$\sup_{g\in G }\sum_{\gamma\in\Gamma} \ee^{Q\dist(\gamma gK,\orig)}=:D<\infty
$.
For $k\in \nat_0$ we have
\begin{eqnarray*}
p_{Y;R,\Omega^k}(\pi_\ast f)&= &\sup_{g\in G}\ee^{R \dist_Y(\Gamma \orig,\Gamma
gK)} |\sum_{\gamma\in \Gamma} f(\gamma g\Omega^kK)|\\
&\le & \sup_{g\in G}\ee^{R \dist_Y(\Gamma \orig,\Gamma gK)}\sum_{\gamma\in
\Gamma}
|f(\gamma  \Omega^kgK)|\\
&\le& \sup_{g\in G}\ee^{R \dist_Y(\Gamma \orig,\Gamma gK)}\sum_{\gamma\in
\Gamma}
\ee^{ (Q-R)\dist(\gamma gK,\orig)} p_{(R-Q),\Omega^k}(f)\\
&\le & \sup_{g\in G}  p_{(R-Q),\Omega^k}(f)
\sum_{\gamma\in \Gamma}
\ee^{ Q \dist(\gamma gK,\orig)}\\
&\le&  D p_{(R-Q),\Omega^k}(f)\ .
\end{eqnarray*}
This estimate shows the continuity of $\pi_\ast$.
Let $\pi^*:S_\infty\cEp_Y\rightarrow S_\infty\cEp$ be the adjoint of $\pi_*$.

To prove the Lemma it suffices to show that
if $f\in  S_\infty\cEp_Y$ with  $F=Bf\in S_\infty\cE_Y$, then $f\in
S_\infty\cE_Y$.
Let $W$ be the parametrix of $B$ constructed in the proof of Lemma \ref{reg}.
Then we have $\pi^\ast f=WF + S\pi^\ast f$ (viewing $F\in {}^\Gamma
S_\infty\cE$).
We have already shown that $WF,S\pi^*f\in S_\infty\cE$.
Since $W,S$ are $G$-equivariant we obtain $\pi^*f\in {}^\Gamma S_\infty\cE $.
This finishes the proof of the Lemma. \hB
\begin{lem}\label{clos}
The operator $B^*:S_{-\infty}\cE_Y\rightarrow S_{-\infty}\cE_Y$
has closed range.
\end{lem}
\proof
Let $\{f_i\}$ be a sequence in $S_{-\infty}\cE_Y$ such that
$B^*f_i=:h_i\rightarrow h\in S_{-\infty}\cE_Y$.
We are to find $f\in S_{-\infty}\cE_Y$ with $B^*f=h$.
If $V_\mu\not=0$, then
$V_\mu=\ker(B^*:S_{-\infty}\cE_Y\rightarrow S_{-\infty}\cE_Y) $,  and
we can project $f_i$ to the $L^2$-orthogonal complement of $V_\mu$.
Thus we can assume that $f_i\perp V_\mu$.

We can assume that for some $R\in \R$ the sequence $p_{Y;R,1}(f_i)$ is bounded.
If not, we divide $f_i$ by $p_{Y;R,1}(f_i)$ obtaining a sequence $\tilde{f}_i$
with  $p_{Y;R,1}(\tilde{f}_i)=1$ and $ B^*\tilde{f}_i\to 0$.
We show below that $\tilde{f}_i$ has a subsequence converging to $F\in
S_{-\infty}\cE_Y$.
Now $ B^*F=0$ and  $F\perp V_\mu$.
Hence $F=0$  contradicting $p_{Y;R,1}(F)=1$.

Since for all  $k\in\nat_0$ the sequence $p_{Y;R,\Omega^k}(h_i)$ is bounded we
conclude
that $p_{Y;R,\Omega^k}(f_i)$ is bounded, too.

Consider a cusp $c_l$ of $Y$ associated to the minimal  parabolic subgroup
$P=MAN\subset G$.
Then a neighbourhood $U$ of infinity of this cusp can be identified with
$\Gamma\cap N\backslash N \times [d_l,\infty)$, $d_l\in A$ large.
If $F\in \cE_Y$, then we define the constant term $F_P\in C^\infty(U,E_{Y|U})$
by
$$F_P(na)=\frac{1}{\vol(\Gamma\cap N\backslash N)}\int_{\Gamma\cap N\backslash
N} F(n^\prime na) dn^\prime \ . $$
Let $\chi\in C^\infty(Y)$ be a cut-off function being one  on $\Gamma\cap
N\backslash N\times [d_l+1,\infty)$ and
zero outside of $U$. Then $F_{c_l}=\chi F_P \in \cE_Y$.
If $Y$ has the cusps $c_l$, $l=1,\dots,r$, then we set $F_c:=\sum_{l=1}^r
F_{c_l}$.

We apply this construction to our sequence $f_i$ obtaining a decomposition
$f_i=f_{i,c}+f_{i,r}$ with $f_{i,r}:=f_i-f_{i,c}$.
Since $f_i$ is bounded in $S_R\cE_Y$, by a Lemma
of Gelfand \cite{franke95}, Thm. 5, the sequence $f_{i,r}$ is
bounded in $S_{-\infty}\cE_Y$.
Since for $R<R^\prime$ the embedding
$S_R\cE_Y\hookrightarrow S_{R^\prime}\cE_Y$ is compact, any bounded
sequence in $S_{-\infty}\cE_Y$ has a converging subsequence.
Thus by taking a subsequence we can assume that $f_{i,r}$ converges in
$S_{-\infty}\cE_Y$.

For
$\dist_Y(\Gamma\orig,y)\ge \max_l d_l+1$ we have $(B^*f)_c(y)= B^*f_c(y)$.
Consider again the cusp $c_l$ and the coordinates $ \Gamma\cap N\backslash
N\times  [d_l,\infty)$.
There are commuting $x,y\in \End(V_\gamma)$ such that
$$( B^*F_{c_l})(n,a)=-((\frac{d}{da}+x)(\frac{d}{da}+y)F_{c_l})(n,a)\ ,$$
$a>d_l+1$.
Assume that $R< - \max( \|x\|) , \|y\| )  $. We set for $F\in S_R\cE_Y$
$$(H_{c_l}F)(n,a):=-\chi(a)  \ee^{-ya}    \int_{a}^\infty
\ee^{(y-x)a_1}\int_{a_1}^\infty \ee^{ xb} F_{c_l}(n,b) db da_1\ .$$
Again, if $Y$ has cusps $c_l$, $l=1,\dots, r$, then we set $H_c:=\sum_{l=1}^r
H_{c_l}$.
Then $H_c:S_R\cE_Y\rightarrow S_R\cE_Y$ is continuous,
and $\supp (H_c B^* F-F_c)\subset V$, where $V\subset Y$ is compact and
independent of $F$.
The proof of continuity is similar to the corresponding argument
in the proof of Lemma \ref{iter}.
Notice that
$\Omega^k H_c-H_c \Omega^k =W$
is a continuous operator $W : S_R\cE_Y\rightarrow C_c^\infty(Y,E_Y)$
and $\supp Wf \subset V^\prime$, where $V^\prime\subset Y$ is compact and
independent of $f$.

Thus $F_i:= H_c B^* h_i-f_{i,c}$ is a bounded sequence in $S_{-\infty}\cE_Y$.
Hence $F_i$ has a subsequence converging in $S_{-\infty}\cE_Y$.
Since $H_c B^* h_i$ converges in $S_{-\infty}\cE_Y$
by taking a subsequence we can assume that $f_{i,c}$ converges
in $S_{-\infty}\cE_Y$, too.
Let $f$ be the limit of $f_i=f_{i,c}+f_{i,r}$ for this subsequence.
Then $ B^*f=h$. This finishes the proof of the lemma.
\hB

We now finish the proof of Theorem \ref{koko}.
By Lemma \ref{clos}  we have
$$\dim \:V_\mu=\dim \ker ( B^*:  S_{-\infty}\cE_Y\rightarrow
S_{-\infty}\cE_Y)=\dim\coker(
B:S_\infty\cEp_Y\rightarrow S_\infty\cEp_Y)\ .$$
By Lemma \ref{smmo}
we have
$$\dim\coker(B:S_\infty\cE_Y\rightarrow S_\infty\cE_Y)\le \dim V_\mu\ .$$
Combining this with Lemma \ref{larger} we obtain the theorem.
\hB

\section{$\Gamma$-cohomology}\label{ggkk}

In this section we discuss properties of the $\Gamma$-cohomology
of distribution vector globalizations of admissible representations
of $G$.

Let $(\pi,V_{\pi,K})\in \hc$ and $\Gamma\subset G$ be a discrete torsion-free
subgroup
of finite covolume.
Let $B=(\Omega-\mu)^k$ such that $BV_{\pi,K}=0$.
\begin{prop}\label{ug1}
We have
$$\dim H^p(\Gamma,V_{\pi,-\infty})<\infty,\quad \forall p\ge 0.$$
\end{prop}
\proof
Let
\begin{equation}\label{ss.ss}
0\rightarrow V_{\pi,-\infty}\rightarrow S_\infty\cE_0\stackrel{\scriptsize
\left(\begin{array}{c}D_0\\B\end{array}\right)}
{\longrightarrow}\begin{array}{c}S_\infty\cE_1\\ \oplus\\
S_\infty\cE_0\end{array}
                                          \stackrel{\scriptsize
\left(\begin{array}{cc}D_1&H_0\\-B&D_0\end{array}\right)}
{\longrightarrow}\begin{array}{c}S_\infty\cE_2\\ \oplus\\
S_\infty\cE_1\end{array}
                                           \stackrel{
\scriptsize\left(\begin{array}{cc}D_2&H_1\\B&D_1\end{array}\right)}
{\longrightarrow}\dots
\end{equation}
be a standard resolution (see Proposition \ref{stare}) of $V_{\pi,-\infty}$.
It is a $\Gamma$-acyclic resolution of $V_{\pi,-\infty}$ by Theorem \ref{gaz}.
The cohomology of the subcomplex of $\Gamma$-invariant vectors is isomorphic to
$H^\ast(\Gamma,V_{\pi,-\infty})$.

For any locally homogeneous vector  bundle $E_Y\rightarrow Y$
we denote by $\cE_Y(B)_{cusp}$ the space of cusp-forms in $\cE_Y(B)$
and by $S_\infty\cE_Y(B)$ the kernel of  $B$ in $S_\infty\cE_Y$.
We consider the subcomplex of the complex of $\Gamma$-invariants of
(\ref{ss.ss})
\begin{equation}\label{rr.rr}
0  \rightarrow S_\infty\cE_{0,Y}(B) \stackrel{\scriptsize
\left(\begin{array}{c}D_0\\0\end{array}\right)}
{\longrightarrow}\begin{array}{c}S_\infty\cE_{1,Y}(B) \\ \oplus\\
\cE_{0,Y}(B)_{cusp}\end{array}
                                          \stackrel{\scriptsize
\left(\begin{array}{cc}D_1&H_0\\0&D_0\end{array}\right)}
{\longrightarrow}\begin{array}{c}S_\infty\cE_{2,Y}(B) \\ \oplus\\
\cE_{1,Y}(B)_{cusp}\end{array}
                                           \stackrel{
\scriptsize\left(\begin{array}{cc}D_2&H_1\\0&D_1\end{array}\right)}
{\longrightarrow}\dots
\end{equation}
and claim that its cohomology is $H^\ast(\Gamma,V_{\pi,-\infty})$.
In fact, let $(f_i,f_{i-1})\in S_\infty\cE_{i,Y}\oplus S_\infty\cE_{i-1,Y}$
be a cochain in the complex of $\Gamma$-invariants of (\ref{ss.ss}).
By the results of Section \ref{klop} there is a unique decomposition
$f_{i-1}=f_{i-1}^{im}+f_{i-1}^{cusp}$, where $f_{i-1}^{im}\in B
S_\infty\cE_{i,Y}$
and $f_{i-1}^{cusp}\in \cE_{i-1,Y}(B)_{cusp}$.
Notice that $\coker B = \coker (\Omega-\mu)$.
Thus modulo a coboundary the cochain $(f_i,f_{i-1})$
is equivalent to $(\tilde{f_i},f_{i-1}^{cusp})$. If in addition $(f_i,f_{i-1})$
and hence $(\tilde{f_i},f_{i-1}^{cusp})$ is a cocycle, then
$B\tilde{f_i}=0$, since $D_{i-1}:\cE_{i-1,Y}(B)_{cusp}\rightarrow
\cE_{i,Y}(B)_{cusp}$
and the range of $B:S_\infty\cE_{i,Y}\rightarrow S_\infty\cE_{i,Y}$ is
transverse
to $\cE_{i,Y}(B)_{cusp}$. We conclude that the cohomology of (\ref{rr.rr})
surjects
onto $H^*(\Gamma,V_{\pi,-\infty})$.

Assume now that $f_i=H_{i-2}g_{i-2}+D_{i-1}g_{i-1}$,
$f_{i-1}=(-1)^{i-1} Bg_{i-1}+D_{i-2}g_{i-2}$, and
$f_{i-1}\in \cE_{i-1,Y}(B)_{cusp}$, $g_{i-2}\in \cE_{i-2,Y}(B)_{cusp}$.
It follows that $Bg_{i-1}=0$ and thus $(g_{i-1},g_{i-2})$ is
a $i-1$-cochain of (\ref{rr.rr}). This show that the cohomology of
(\ref{rr.rr})
maps injectively to $H^*(\Gamma,V_{\pi,-\infty})$ proving the claim.

The lemma now follows since (\ref{rr.rr}) is a complex of finite-dimensional
vector spaces (see e.g. \cite{harishchandra68}, Thm.1).
\hB
If $\Gamma$ is cocompact, then we can prove a Poincar\'e duality theorem as in
\cite{bunkeolbrich947}, Proposition 5.2.
It relates the $\Gamma$-cohomology of the distribution vector globalization
of an admissible representation of $G$ with the $\Gamma$-cohomology
of the smooth globalization of its dual.
\begin{prop} \label{ug2}
Let $\Gamma$ be cocompact. The $\Gamma$-cohomology of $V_{\pi,\pm\infty}$
satisfies Poincar\'e duality
$$H^p(\Gamma,V_{\pi,-\infty})^\ast\cong H^{n-p}(\Gamma,V_{\tilde\pi,\infty}),$$
where $n=\dim(X)$.
\end{prop}
Let $V_{\pi,\pm\omega}$ be the minimal and maximal globalizations of
$V_{\pi,K}$,
respectively. If $\Gamma$ is cocompact, then
$\cE_Y(B)_{cusp}=\cE_Y(B)=S_\infty\cE_Y(B)$.
Combining \cite{bunkeolbrich947}, Proposition 5.1 with Proposition \ref{ug1} we
obtain
\begin{kor}
Let $\Gamma$ be cocompact. Then for $(\pi,V_{\pi,K})\in\hc$ we have
\begin{eqnarray*}
H^\ast(\Gamma,V_{\pi,\omega})&=&H^\ast(\Gamma,V_{\pi,\infty})\\
H^\ast(\Gamma,V_{\pi,-\infty})&=&H^\ast(\Gamma,V_{\pi,-\omega})\ ,
\end{eqnarray*}
where the identifications are induced by the natural inclusions
of the globalizations.
\end{kor}
If $\Gamma$ has torsion, then there exists a cofinite torsion free  normal
subgroup
$\Gamma^\prime\subset \Gamma$. Employing the spectral sequence
for the group cohomology associated to the extension
$$0\rightarrow \Gamma^\prime\rightarrow \Gamma\rightarrow F\rightarrow 0\ ,$$
where $F$ is some finite group, one easily obtains   a generalization of the
results
of the present section to $\Gamma$.
The spectral sequence degenerates at the second term  since
the higher cohomology of a finite group with
coefficients in a vector spaces over  $\C$ vanishes.
It follows $H^*(\Gamma,V)=H^*(\Gamma^\prime,V)^F$ for any $\Gamma$-module
$V$ over $\C$.
 \section{Fuchsian groups of the first kind}\label{fufu}

In this section we give a detailed discussion of the $\Gamma$-cohomology
 of a Fuchsian group of the first kind with coefficients
in the distribution vector globalization of principal series
representations. In this case we know explicit standard resolutions.

Let $\Gamma\subset PSL(2,\R)=:G$ be  a discrete torsion-free subgroup
of finite covolume. Such a $\Gamma$ is called a Fuchsian group of the first
kind
and it acts freely on the hyperbolic plane $X=H^2$.
The quotient $Y=\Gamma\backslash X$ is a complete Riemann surface of finite
volume.

The group $G$ acts on the circle $S^1$
which can be identified with the
boundary $\partial X$ of $X$ using the Poincar\'e disc model.
Let $T\rightarrow S^1$ be the complexified tangent bundle of $S^1$.
It is $G$-homogeneous and we can form complex powers $T^\lambda\rightarrow
S^1$, $\lambda\in\C$.
The  number $\lambda\in\C$ parametrizes a  principal series representation
$(\pi^\lambda,H^\lambda)$ of $G$ on the Hilbert space
$L^2(S^1,T^{\lambda-1/2})$.
By $H^\lambda_{-\infty}$ we denote the space of its distribution vectors.

For $\lambda\not=-1/2,-3/2,-5/2,\dots$
combining a theorem of Helgason
(\cite{helgason70}, \cite{helgason84} Introduction Thm. 4.3)
with the characterization of the distribution vector globalization
(Wallach \cite{wallach92} Ch. 11, Casselman \cite{casselman89}) we see that
the Poisson transform $P_\lambda$ is an $G$-equivariant isomorphism
$$P_\lambda:H^\lambda_{-\infty}\stackrel{\sim}{\rightarrow} S_\infty\cE(B)\ ,$$
where $B=\Omega-1/4+\lambda^2$ and $E=X\times\C$ is the trivial bundle
(in this special situation this fact was first obtained by Lewis
\cite{lewis78}).
Thus a standard resolution of the principal series representation
$H^\lambda_{-\infty}$
for $\lambda\not=-1/2,-3/2,-5/2,\dots$
is simply
$$0\rightarrow H^\lambda_{-\infty}\stackrel{P_\lambda}{\rightarrow}S_\infty\cE
\stackrel{B}{\rightarrow}S_\infty\cE\rightarrow 0\ .$$
The complex (\ref{rr.rr}) reduces to
\begin{equation}\label{red}0\rightarrow
S_\infty\cE_Y(B)\stackrel{0}{\rightarrow}\cE_Y(B)_{cusp}\rightarrow 0\
.\end{equation}
\begin{prop}\label{gener}
For $\lambda\not=-1/2,-3/2,-5/2,\dots$ we have
\begin{eqnarray*}
H^0(\Gamma,H^\lambda_{-\infty})&=&S_\infty\cE_Y(B)\\
H^1(\Gamma,H^\lambda_{-\infty })&=&\cE_Y(B)_{cusp}\\
H^2(\Gamma,H^\lambda_{-\infty })&=& 0.
\end{eqnarray*}
Moreover, $\chi(\Gamma,H^\lambda_{-\infty})=\dim
H^0(\Gamma,H^\lambda_{-\infty})-\dim H^1(\Gamma,H^\lambda_{-\infty })=r$,
where $r$ is the number of cusps of $Y$.
If $H^1(\Gamma,H^\lambda_{-\infty})\not=0$, then $\lambda\in
\imath\R\cup(-1/2,1/2)$.
\end{prop}
\proof
The first part of the Proposition follows immediately from (\ref{red}).
$\cE_Y(B)_{cusp}\subset \ker_{L^2}(B:\cE_Y\rightarrow \cE_Y)$ and
$\spec_{L^2}\Delta_Y\subset [0,\infty)$
implies the last assertion.
Let $p:S_\infty\cE_Y(B)\rightarrow \C^{2r}$ be the linear map taking the
constant
term. For any cusp the constant term has two components
(the incoming and the outcoming).
It is known that $\dim\im(p)=r$.
In fact, the range of $p$ is generated by the constant terms of regular
Eisenstein
series and their residues. The scattering matrix fixes the relation
between the two components of the constant term.
Since $\ker(p)=\cE_Y(B)_{cusp}$, the assertion about the Euler characteristic
 follows. \hB

Now we discuss the case $\lambda=-k/2$, $k=1,3,5\dots$
taking the structure of $H^{\pm k/2}_{-\infty}$ as a $G$-module into account.
We first exploit the exact  sequence
\begin{equation}\label{wert}0\rightarrow F_k \rightarrow
H^{k/2}_{-\infty}\rightarrow D^+_k\oplus D^-_k\rightarrow 0\ ,\end{equation}
where $F_k$ is the finite-dimensional representation of $G$ of dimension $k$
and $D^\pm_k$ are  the distribution vectors  of
holomorphic and anti-holomorphic discrete series representations.
Let $r$ denote the number of cusps of $Y$ and $g$ denote the genus.
We have $H^i(\Gamma,H^{k/2}_{-\infty})=0$, $i\ge 1$, and
$\dim H^0(\Gamma,H^{k/2}_{-\infty})=\dim S_\infty\cE_Y(B) = r$ by Proposition
\ref{gener}.

The long exact cohomology sequence associated to  (\ref{wert}) gives
\begin{eqnarray*}
&&0\rightarrow H^0(\Gamma,F_k)\rightarrow
H^0(\Gamma,H^{k/2}_{-\infty})\rightarrow H^0(\Gamma,D_k^+\oplus D_k^-)
\stackrel{\delta}{\rightarrow} H^1(\Gamma,F_k)\rightarrow 0\\
&&\hspace{2cm}H^1(\Gamma,D_k^+\oplus D_k^-) = 0   \ .
\end{eqnarray*}
An investigation of the long exact sequence associated to
$$0\rightarrow D_k^+\oplus D_k^-\rightarrow H^{-k/2}_{-\infty}\rightarrow
F_k\rightarrow 0$$
leads to
 \begin{eqnarray*}&&0\rightarrow H^0(\Gamma,D_k^+\oplus D_k^-)\rightarrow
H^0(\Gamma,H^{-k/2}_{-\infty}) \rightarrow  H^0(\Gamma,F_k)\rightarrow 0\\
&&\hspace{2cm}H^1(\Gamma,H^{-k/2}_{-\infty})=H^1(\Gamma,F_k)\ .\end{eqnarray*}
  \begin{prop}\label{wwqqq}
We have
\begin{eqnarray*}
\dim H^0(\Gamma,D_k^+\oplus D_k^-)&=&r+\dim H^1(\Gamma,F_k)-\dim
H^0(\Gamma,F_k)\\
\dim H^0(\Gamma,H^{-k/2}_{-\infty})&=&r+ \dim H^1(\Gamma,F_k)\\
\dim H^1(\Gamma,H^{-k/2}_{-\infty})&=&\dim H^1(\Gamma,F_k)\\
\chi(\Gamma,H^{-k/2}_{-\infty})&=&r
\end{eqnarray*}
\end{prop}
Since $r>1$ or $g>0$ we have for $k>1$
\begin{eqnarray*}
H^0(\Gamma,F_1)&=&1\\
H^1(\Gamma,F_1)&=&2g+r-1\\
H^0(\Gamma,F_k)&=&0\\
H^1(\Gamma,F_k)&=&k(2g+r-2)\ .
\end{eqnarray*}
It follows that
\begin{eqnarray*}
\dim H^0(\Gamma,D_1^+\oplus D_1^-)&=&2g+2r-2 \\
\dim H^0(\Gamma,D_k^+\oplus D_k^-)&=&k(2g -2)+(k+1)r\\
\dim H^0(\Gamma,H^{-1/2}_{-\infty})&=&2g+2r-1\\
\dim H^1(\Gamma,H^{-1/2}_{-\infty})&=&2g+r-1\\
\dim H^0(\Gamma,H^{-k/2}_{-\infty})&=&k(2g-2)+(k+1)r\\
\dim H^1(\Gamma,H^{-k/2}_{-\infty})&=&k(2g+r-2) \ .
\end{eqnarray*}

The aim of the following discussion is to understand Proposition \ref{wwqqq}
in the framework of standard resolutions.

Let $\K$ be the canonical bundle of $X$ (viewing $X$ as a complex  manifold)
and $\K^i$ be its $i$'th power.
\begin{lem}\label{reso}
A standard resolution of $H^{-k/2}_{-\infty}$ is given by
\begin{eqnarray*}
\lefteqn{0\longrightarrow H^{-k/2}_{-\infty}\stackrel{P}{\longrightarrow}
S_\infty\cK^{(k+1)/2}\oplus S_\infty\cK^{-(k+1)/2}}\hspace{3cm}\\
&&\stackrel{\scriptsize\left(\begin{array}{c}
\bar{\partial}^{(k+1)/2}\oplus\partial^{(k+1)/2}\\\Omega+
\frac{k^2-1}{4}\end{array}\right)}{
\longrightarrow}\begin{array}{c}S_\infty\cK^0\\\oplus\\
 S_\infty\cK^{(k+1)/2}\oplus
S_\infty\cK^{-(k+1)/2} \end{array}
\longrightarrow\\
&&\stackrel{\scriptsize
(-\Omega-\frac{k^2-1}{4},\bar{\partial}^{(k+1)/2}
\oplus\partial^{(k+1)/2})}{\longrightarrow}
S_\infty\cK^0\longrightarrow 0\ .
\end{eqnarray*}
 Here $\bar{\partial}:\cK^{(l+1)/2}\rightarrow \cK^{(l-1)/2}$, $l>0$ odd,
is the contraction of the anti-holomorphic part of the canonical connection
with the K\"ahler form and similarly
$\partial:\cK^{-(l+1)/2}\rightarrow\cK^{-(l-1)/2}$
is the contraction of the holomorphic part of the canonical connection with the
K\"ahler form.
In abuse of notation we write $\bar{\partial}^{(k+1)/2}, \partial^{(k+1)/2}$
for the composition of the corresponding number of $\bar{\partial}$'s,
$\partial$'s,
respectively.
 \end{lem}
\proof
Define $B:=\Omega+\frac{k^2-1}{4}$.
As $G$-modules the eigenspaces
$S_\infty\cK^{\pm(k+1)/2}(B)$ have three-step composition  series'
with composition factors $D^\pm_k,F_k,D^\mp_k$ (see \cite{olbrichdiss}).
The operators $\bar{\partial}^{(k+1)/2}$, $\partial^{(k+1)/2}$, respectively,
annihilate exactly the first composition factor.
Above we have seen that
$S_\infty\cK^0(B) = H^{k/2}_\infty$ with composition factors $D_k^+\oplus
D_k^-,F_k$.
This discussion   shows that the complex above resolves
$H^{-k/2}_\infty$ on a ${\bf K}$-theoretic level.
In order to show that the complex is in fact a resolution of the specific
extension
$H^{-k/2}_\infty$  of the composition factors note that there is an injective
Poisson transform $P$.
\hB
Now we take the $\Gamma$-invariant vectors in the standard resolution.
The complex (\ref{rr.rr}) reduces to
\begin{eqnarray*}
\lefteqn{0\longrightarrow S_\infty\cK^{(k+1)/2}_Y(B)\oplus
S_\infty\cK^{-(k+1)/2}_Y(B)}\hspace{3cm}\\
&&\stackrel{\scriptsize\left(\begin{array}{c}
\bar{\partial}^{(k+1)/2}\oplus\partial^{(k+1)/2}\\0
\end{array}\right)}{\longrightarrow}\begin{array}{c}
S_\infty\cK_Y^0(B)\\\oplus\\ \cK^{(k+1)/2}_Y(B)_{cusp}
\oplus \cK^{-(k+1)/2}_Y(B)_{cusp} \end{array}
\longrightarrow\\
&&\stackrel{\scriptsize(0,\bar{\partial}^{(k+1)/2}
\oplus\partial^{(k+1)/2})}{\longrightarrow} \cK^0_Y(B)_{cusp}
\longrightarrow 0\ .
\end{eqnarray*}
In order to make further reductions we assume that $Y$ is not compact
(the case of compact $Y$ is an easy exercise and left to the reader).
Since for $k>1$ the operator $B$ on $L^2(Y)$ is strictly positive    and
$$\cK^0_Y(B)_{cusp}\subset \ker_{L^2}(B:\cK^0_Y\rightarrow \cK^0_Y)$$ for $k>1$
we have $\cK^0_Y(B)_{cusp}=0$. This is also true for $k=1$ since then
$\ker_{L^2}(B:\cK^0_Y\rightarrow \cK^0_Y)$ is generated  by the constant
function which
is not a cusp form.

We claim that for $k>1$
\begin{eqnarray}0&=&\bar{\partial}^{(k+1)/2}:
S_\infty\cK^{(k+1)/2}_Y(B)\rightarrow S_\infty
\cK^{0}_Y(B)\label{firstl}\\
0&=&\partial^{(k+1)/2}:S_\infty\cK^{-(k+1)/2}_Y(B)
\rightarrow S_\infty\cK^0_Y(B)\ .\label{secondl}
\end{eqnarray}
To prove (\ref{firstl}) we argue that already
$$ 0=\bar{\partial}:S_\infty\cK^{(k+1)/2}_Y(B)\rightarrow
S_\infty\cK^{(k-1)/2}_Y(B)\ .$$
In upper half-plane coordinates $(x,y)$, $y>0$, $z=x+\imath y$,
we trivialize $\K^{(k+1)/2}$ using the section $(\frac{dz}{y})^{(k+1)/2}$.
Then the operator $B$ has the form $y^2\Delta + \imath  (k+1)
\frac{\partial}{\partial x}+\frac{k^2-1}{4}$.
If $\phi\in S_\infty\cK^{(k+1)/2}(B)$, then its constant term is a linear
combination of
$y^{-k}dz^{(k+1)/2}$ and $dz^{(k+1)/2}$. It follows that $\bar{\partial}\phi\in
L^2$.
But $\ker_{L^2}(B)=0$ on $\K^{ (k-1)/2}$, $k >1$.
Equation (\ref{secondl}) follows
by complex conjugation.

 The spaces
\begin{eqnarray*}\cA_{k+1}&:=&\ker(\bar{\partial}:
S_\infty\cK^{(k+1)/2}_Y(B)\rightarrow S_\infty\cK^{(k-1)/2}_Y(B))\\
\bar{\cA}_{k+1}&:=&\ker( \partial: S_\infty\cK^{-(k+1)/2}_Y(B)\rightarrow
S_\infty\cK^{-(k-1)/2}_Y(B))
\end{eqnarray*}
 are called the holomorphic and anti-holomorphic automorphic forms of weight
$k+1$.
The subspaces
\begin{eqnarray*}
\cS_{k+1}&:=&\cK^{(k+1)/2}_Y(B)_{cusp}\subset \cA_{k+1}\\
\bar{\cS}_{k+1}&:=&\cK^{-(k+1)/2}_Y(B)_{cusp}\subset \bar{\cA}_{k+1}
\end{eqnarray*}
are the holomorphic and anti-holomorphic cusp-forms of weight $k+1$.
Since for $k>1$  the space $S_\infty\cK^0_Y(B)$ is generated by Eisenstein
series
its dimension is equal to the number of cusps $r$.
For $k=3,5,\dots$ we have
\begin{eqnarray*}
\dim H^0(\Gamma,H^{-k/2}_{-\infty})&=&\dim(\cA_{k+1} \oplus \bar{\cA}_{k+1})\\
\dim H^1(\Gamma,H^{-k/2}_{-\infty})&=&\dim(\cS_{k+1} \oplus \bar{\cS}_{k+1})+r\
{}.
\end{eqnarray*}
Employing $\dim \cA_{k+1} = \dim \bar{\cA}_{k+1}  = k(g-1) +r (k+1)/2$
and $\dim \cS_{k+1} = \dim \bar{\cS}_{k+1}  = k(g-1) + r(k-1)/2$ (see
\cite{shimura94}, Thm. 2.23) we recover the result of  Proposition \ref{wwqqq}.

We now consider the case $k=1$. In this case both
\begin{eqnarray}&&\bar{\partial}   :S_\infty\cK^{1}_Y(B)\rightarrow
S_\infty\cK^{0}_Y(B)\ , \\
&&\partial :S_\infty\cK^{-1}_Y(B)\rightarrow S_\infty\cK^0_Y(B)
\end{eqnarray}
 have a one-dimensional range spanned by the constant function.
Note  the $r$-dimensional space $ S_\infty\cK^0_Y(B)$ is generated by the
regular Eisenstein series and the constant function.
We obtain for $k=1$
\begin{eqnarray*}
\dim H^0(\Gamma,H^{-1/2}_{-\infty})&=&\dim(\cA_2 \oplus \bar{\cA}_2)+1\\
\dim H^1(\Gamma,H^{-1/2}_{-\infty})&=&\dim(\cS_2 \oplus \bar{\cS}_2)+r-1\ .
\end{eqnarray*}
Using $\dim \cA_2=\dim \bar{\cA}_2=g-1+r$, $ \dim \cS_1 = \dim \bar{\cS}_2  =
g$ (see \cite{shimura94}, Thm. 2.23)
we again recover the results of Proposition \ref{wwqqq}.

We finish this section with a discussion of the relation of the
sequence (\ref{wert}) with the Eichler homomorphism
 $E:\cS_{k+1} \rightarrow H^1(\Gamma,F_k)$
(see e.g. \cite{shimura94}).
We first recall its definition.
If we restrict $F_k$ to the maximal compact subgroup $K=S^1\subset G$,
then it decomposes as $(F_k)_{|K}=\sum_{m=1}^k \C_{ m-(k+1)/2} $, where
$\C_l$ is the representation of $S^1$ on $\C$ given by $z\mapsto z^l$.
Note that $\K^l=G\times_K \C_{-l}$. The map
$(gK,f)\mapsto [g,g^{-1}f]$, $g\in G$, $f\in F_k$ defines an isomorphism
$$G/K\times F_k\stackrel{\sim}{\rightarrow} G\times_K (F_k)_{|K}=\oplus_{m=1}^k
\K^{ m-(k+1)/2}\ .$$
We denote the canonical $G$-equivariant embedding
$\K^{\pm(k-1)/2}\hookrightarrow X\times F_k\cong L^2\otimes F_k$
by $j_\pm$, where the second identification is given by the Hodge-$*$-operator
and $L^l$ is the bundle of $l$-forms on $X$. Taking the tensor product with
$\K^{\pm 1}$ we obtain embeddings
$i_\pm:\K^{\pm(k+1)/2}\hookrightarrow L^1\otimes F_k$.
If $f\in \cS_{k+1}=\cK^{(k+1)/2}(B)_{cusp}$, then $i_+(f)$ is a closed
form and represents $E(f)$.

The long exact cohomology sequence associated to (\ref{wert}) induces a
boundary map
$$\delta:H^0(\Gamma,D^+_k\oplus D^-_k)\rightarrow H^1(\Gamma,F_k)\ .$$

\begin{lem}\label{eichler}
The map $\delta$ restricted to the cusp forms coincides with $E$.
\end{lem}
\proof
In order to compute $\delta$ we need a $\Gamma$-acyclic resolution of
(\ref{wert})
as a sequence.
One possibility is the following:
\begin{equation}{\scriptsize \label{saur}
\begin{array}{ccccccccc}
&&0&&0&&&&\\
&&\uparrow&&\uparrow&&&&\\
0&\rightarrow&S_\infty\cL^2\otimes
F_k&\stackrel{=}{\rightarrow}&S_\infty\cL^2\otimes F_K&\rightarrow&0&&\\
&&d\uparrow&&a\uparrow&&\uparrow&&\\
0&\rightarrow&S_\infty\cL^1\otimes F_k&\rightarrow&S_\infty\cL^1\otimes
F_K\oplus S_\infty\cK^{(k-1)/2}\oplus S_\infty\cK^{-(k-1)/2}&\rightarrow
&S_\infty\cK^{(k-1)/2}\oplus S_\infty\cK^{-(k-1)/2}&\rightarrow&0\\
&&d\uparrow&&b\uparrow& &c\uparrow&&\\
0&\rightarrow&S_\infty\cL^0\otimes F_k&\rightarrow&S_\infty\cL^0\otimes
F_K\oplus S_\infty\cK^{(k+1)/2}\oplus S_\infty\cK^{-(k+1)/2}&\rightarrow
&S_\infty\cK^{(k+1)/2}\oplus S_\infty\cK^{-(k+1)/2}&\rightarrow&0\\
&&\uparrow&&\uparrow& &\uparrow&&\\
0&\rightarrow& F_k&\rightarrow &H^{k/2}&\rightarrow &D^+_k\oplus D^-_k
&\rightarrow&0\\
&&\uparrow&&\uparrow&&\uparrow&&\\
&&0&&0&&0&&\ .
\end{array}}
\end{equation}
Here the maps are given by
\begin{eqnarray*}
a&:=&d-j_+-j_-\\
b&:=&\left(\begin{array}{ccc}d&i_+&i_-
\\0&\bar{\partial}&0\\0&0&\partial\end{array}\right)\\
c&:=&\bar{\partial}\oplus\partial\ .
 \end{eqnarray*}
The horizontal maps are the obvious embeddings and projections, respectively.
We leave the verification of the commutativity of the diagram and the exactness
of the middle column to  the reader.

The boundary map $\delta$ is now obtained by the usual diagram chasing. Let
$[\alpha]\in H^0(\Gamma,D^+_k)$ be represented by a holomorphic $\Gamma$-
invariant form $\alpha\in{}^\Gamma S_\infty \cK^{(k+1)/2}$. Then $\delta
[\alpha] \in H^1(\Gamma,F_k)$
 is represented by the closed form $ i_+(\alpha)\in S_\infty \cL^1\otimes F_k$.
Comparing this with the definition of the Eichler map $E$ given above
we finish the proof of the lemma.
\hB

\bibliographystyle{plain}

\end{document}